\ifcase2
    \documentclass[a4paper,twocolumn]{article}
    \usepackage[a4paper,landscape]{geometry}
\or \documentclass[12pt,a5paper]{article}
\or \documentclass[12pt,a4paper]{article}
\fi

\usepackage{euler,amsfonts,amsmath,natbib,defns}
\IfFileExists{ajr.sty}{\usepackage{ajr}}
{\usepackage{graphicx,url}\usepackage[dvips,usenames]{color}}
\allowdisplaybreaks

\title{Normal form transforms separate slow and fast modes in
stochastic dynamical systems}

\author{A.~J. Roberts\thanks{Computational Engineering and Sciences
Research Centre, Department of Mathematics \& Computing, University of
Southern Queensland, Toowoomba, Queensland~4352, \textsc{Australia}.
\url{mailto:aroberts@usq.edu.au}}}

\newcommand{\Z}[1]{e^{#1t}\star}
\newcommand{\sgn}{\operatorname{sgn}}
\newcommand{\res}{\operatorname{Res}}
\renewcommand{\P}[1]{Principle~\ref{p:#1}}
\newcommand{\Oeps}[1]{\Ord{\epsilon^{#1}}}
\newcommand{\freq}{\Omega}
\newcommand{\ff}{\omega}
\newcommand{\xint}{\int_D}
\newcommand{\pdf}{\textsc{pdf}}
\newcommand{\ssm}{\textsc{ssm}}
\newcommand{\E}[1]{\operatorname{E}\big[#1\big]}
\newcommand{\cc}[1]{{#1}^*}
\newcommand{\avx}{{\bar x}}
\newcommand{\spde}{\textsc{spde}}
\let\alt\breve
\let\aalt\acute

\begin{document}

\maketitle

\begin{abstract}  
Modelling stochastic systems has many important applications.
Normal form coordinate transforms are a powerful way to untangle
interesting long term macroscale dynamics from detailed microscale
dynamics.
We explore such coordinate transforms of stochastic differential systems when the
dynamics has both slow modes and quickly decaying modes.
The thrust is to derive normal forms useful for macroscopic modelling
of complex stochastic microscopic systems.
Thus we not only must reduce the dimensionality of the dynamics, but
also endeavour to separate all slow processes from all fast time processes, both deterministic and stochastic.
Quadratic stochastic effects in the fast modes contribute to the drift of the important slow
modes. 
The results will help us accurately model, interpret and simulate multiscale stochastic systems.
\end{abstract}

\tableofcontents

\section{Introduction}

Normal form coordinate transformations provide a sound basis for simplifying multiscale nonlinear dynamics \cite[e.g.]{Elphick87b, Cox93b}.  
In systems with fast and slow dynamics, a coordinate transform is sought that decouples the slow from the fast.  
The decoupled slow modes then provide accurate predictions for the long term dynamics.      
Arguably, such normal form coordinate transformations that decouple slow and fast modes provide a much more insightful view of simplifying dynamics than other, more popular, techniques.  
Averaging is perhaps the most popular technique for simplifying dynamics \cite[Chapters 11--13, e.g.]{Verhulst05}, especially for stochastic dynamics that we explore here~\cite[e.g.]{Pavliotis06a, Givon06}.  
But averaging fails in many cases.  
For example, consider the simple, linear, slow-fast, system of stochastic differential equations (\sde{}s)
\begin{equation}
    dx=\epsilon y\,dt \qtq{and}
    dy=-y\,dt+dW\,,
    \label{eq:sde0}
\end{equation}
where for small parameter~$\epsilon$ the variable~$x(t)$ evolves slowly compared to the fast variable~$y(t)$.  
Let us compare the predictions of averaging and a `normal form' coordinate transform.  First consider averaging: the fast variable~$y$, being an Ornstein--Uhlenbeck process, rapidly approaches its limiting \pdf\ that is symmetric in~$y$.  
Then averaging the $x$~equation leads to the prediction $d\bar x=0\,dt$\,; that is, averaging predicts nothing happens.  
Yet the slow $x$~variable must fluctuate through its forcing by the fast~$y$.  
Second, and similar to illuminating coordinate transforms used in Sections~\ref{sec:eg}--\ref{sec:nfsmsde}, introduce new coordinates $X$~and~$Y$ to replace $x$~and~$y$ where
\begin{equation}
    x=X-\epsilon Y+\epsilon \int_{-\infty}^t e^{\tau-t} \,dW_\tau \qtq{and}
    y=Y+\int_{-\infty}^t e^{\tau-t} \,dW_\tau \,.
    \label{eq:sde0nf}
\end{equation}
In the $X$~and~$Y$ coordinates the \sde\ system~\eqref{eq:sde0} decouples to simply
\begin{equation}
    dX=\epsilon\,dW \qtq{and}
    dY=-Y\,dt\,.
    \label{eq:Sde0}
\end{equation}
In these new coordinates $Y\to0$ exponentially fast.  
Thus in the long term the only significant dynamics occurs in the new slow variable~$X$ which system~\eqref{eq:Sde0} shows undergoes a random walk.  
The method of averaging completely misses this random walk: true, the mean~$\bar x$ remains at zero; but the growing spread about the mean is missed by averaging.  Stochastic coordinate transforms such as~\eqref{eq:sde0nf} decouple fast and slow variables to empower us to extract accurate models for the new slow variable~$X$.  
They are called  `normal form' transformations because this decoupling of stochastic dynamics is analogous to corresponding simplifications in deterministic systems \cite[e.g.]{Murdock03, Arnold03}.  
This article establishes useful properties for such stochastic normal form coordinate transformations in modelling multiscale nonlinear stochastic dynamical systems.

Stochastic \ode{}s and \pde{}s have many important applications.
Here we restrict attention to nonlinear \sde{}s when the dynamics of
the \sde\ has both long lasting slow modes and decaying fast
modes~\cite[e.g.]{Arnold98}.
The aim underlying all the exploration in this article is to derive
normal forms useful for macroscopic modelling of stochastic systems
when the systems are specified at a detailed microscopic level.
Thus we endeavour to separate \emph{all} fast time processes from the
slow processes~\cite[e.g.]{Chao95, Roberts05c}.
Such separation is especially interesting in stochastic systems as white noise has fluctuations on \emph{all} time scales.
In contrast, almost all previous approaches have been content to derive
normal forms that support reducing the dimensionality of the dynamics.
Here we go further than other researchers and \emph{additionally separate
fast time processes from the slow modes.}

\cite{Arnold98} developed rigorous theory for stochastic coordinate
transforms to a normal form in the ``non-resonant'' cases.
Then a stochastic system is rigorously equivalent to a decoupled
linearised system no matter how disparate the time scales.
However, stochastic fast-slow systems are essentially resonant.
\cite{Arnold98} report,~\S4, that their results do apply, but that the
resulting normal form is generically nonlinear \cite[\S8.4
also]{Arnold03}.
They comment that the normal form transformation involves
anticipating the noise processes, that is, involving integrals of the
noise over a fast time scale of the future.
However, in contrast to the examples of \cite{Arnold98}
\cite[corrected]{Arnold03}, Sections \ref{sec:eg}~and~\ref{sec:nfsmsde}
argue that such anticipation can always be removed from the slow modes
with the result that no anticipation is required after the fast transients
decay.
Furthermore, Sections \ref{sec:eg}~and~\ref{sec:nfsmsde} argue that on
the stochastic slow manifold all noise integrals can be removed from
terms linear in the noise to leave a slow mode system, such as the
simple $dX=\epsilon\,dW$ of the normal form~\eqref{eq:Sde0}, in which
there are no fast time integrals at all.
The arguments demonstrate that, except for some effects nonlinear in
the noise, all fast time processes can be removed from the slow modes
of a normal form of stochastic systems.

The theory of \cite{Arnold98} applies only to finite dimensional
stochastic systems.
Similarly, \cite{Du06}'s theory of invariant manifold reduction for
stochastic dynamical systems also only applies in finite dimensions.
But many applications are infinite dimensional; for example, the
discretisation of stochastic \pde{}s approximates an inertial manifold
of stochastic dynamics \cite[]{Roberts05c}.
Following the wide recognition of the utility of inertial manifolds
\cite[e.g.]{Temam90}, \cite{Bensoussan95} proved the existence of
attractive stochastic inertial manifolds in Hilbert spaces.
The stochastic slow manifolds obtained in Sections
\ref{sec:eg}~and~\ref{sec:nfsmsde} via stochastic normal forms are
examples of such stochastic inertial manifolds, albeit still in finite
dimensions.

To derive a normal form we have to implement a coordinate
transformation that simplifies a stochastic system.
But the term `simplify' means different things to different people
depending upon how they wish to use the `simplified' stochastic system.
Our aim throughout this article is to create stochastic models that may
efficiently simulate the long term dynamics of multiscale stochastic
systems.
This aim is a little different to that of previous researchers and so
the results herein are a little different.
For example, \cite{Coullet83} and \cite{Arnold98} do not avoid fast time integrals because their aim is different.
Principles that we require are the following:
\begin{enumerate}
\item \label{p:secular} Avoid unbounded (secular) terms in the
transformation and the evolution (ensures uniform asymptotic approximations);

\item \label{p:decouple} Decouple all the slow processes from the fast processes (ensures a valid long term model);

\item \label{p:zero} Insist that the stochastic slow manifold is
precisely the transformed fast modes being zero;
    
\item \label{p:ruthless} Ruthlessly eliminate as many as possible of
the terms in the evolution (to simplify at least the algebraic form of the \sde{}s);

\item \label{p:memory} Avoid as far as possible fast time memory
integrals in the evolution (to endeavour to remove all fast time processes from the slow modes).
    
\end{enumerate}
In general we can meet all these principles, although the last two are
easy as they are only phrased as `as far\slash many as
possible': Section~\ref{sec:eg} explores the issues in a particular
example stochastic system; whereas Section~\ref{sec:nfsmsde} presents
general arguments for finite dimensional, nonlinear, stochastic
differential systems.
Sections \ref{sec:trsnfdm}~and~\ref{sec:aa} additional show that two
other alternative, and superficially attractive, principles are not
useful in the context of macroscale modelling of stochastic systems.

\cite{Srinamachchivaya90} and \cite{Srinamachchivaya91} emphasise the
importance of effects quadratic in the stochastic noise ``in order to
capture the stochastic contributions of the stable modes to the drift
terms of the critical modes.'' Sections
\ref{sec:eg}~and~\ref{sec:nfsmsde} also address such important
quadratic effects.
The generic result of this normal form approach is that not all the
memory integrals can be removed from the evolution of the stochastic 
slow variables: some terms quadratic in the noise retain memory integrals.

Section~\ref{sec:imstm} explores the implications of these results for
macroscale simulation of stochastic systems.  The normal form
approach empowers us to address the effect of anticipatory integrals, 
the influence of the noise on averages, especially noise induced drift,
and the failure of averaging to provide a systematic basis for
macroscale simulation.


Lastly, Section~\ref{sec:osc} discusses in detail a normal form of a stochastically
forced Hopf bifurcation.
Not because it is a Hopf bifurcation, but instead because it is a
generic example of stochastic effects in oscillatory dynamics.
The primary issue is how to `average' over both the nonlinear
oscillation and the noise effects to generate a prescriptive
model of the dynamics over much longer time scales.
A complex valued, time dependent, coordinate transform
can, with considerable care, derives a model \sde\ that is valid for
simulating the long term evolution of the stochastic oscillating
dynamics.
A future application could be to the modelling of atmospheric white
noise forcing of oceanic modes: \cite{Pierce01} discusses this
situation from an oceanographer's perspective.

\section{Explore in detail a simple nonlinear stochastic system}
\label{sec:eg}

This section considers the dynamics of one of the most elementary,
nonlinear, multiscale stochastic systems:
\begin{equation}
    dx=-xy\,dt
    \quad\text{and}\quad
    dy=(-y+x^2-2y^2)dt+\sigma\,dW\,.
    \label{eq:toy0}
\end{equation}
The issues raised, and their resolution, in this relatively simple stochastic system are generic as seen in Section~\ref{sec:nfsmsde}.

Throughout this article I adopt the Stratonovich interpretation of
\sde{}s, as does \cite{Arnold98}, so that the usual rules of calculus apply.
To ease asymptotic analysis I also adopt hereafter the notation of
applied physicists and engineers.
Thus I formally explore the \sde\ system~\eqref{eq:toy0} in the
equivalent form of
\begin{equation}
    \dot x=-xy
    \quad\text{and}\quad
    \dot y=-y+x^2-2y^2+\sigma\phi(t)\,,
    \label{eq:toy}
\end{equation}
where overdots denote formal time derivatives and the `white noise'~$\phi(t)$
is the formal time derivative of the Wiener process~$W(t)$.
Both the Stratonovich interpretation and the adoption of this formal notation
empowers the use of computer algebra to handle the multitude of details.

\begin{figure}[tbp]
    \centering
    \begin{tabular}{c@{}c}
        \raisebox{20ex}{$y$} & 
        \includegraphics{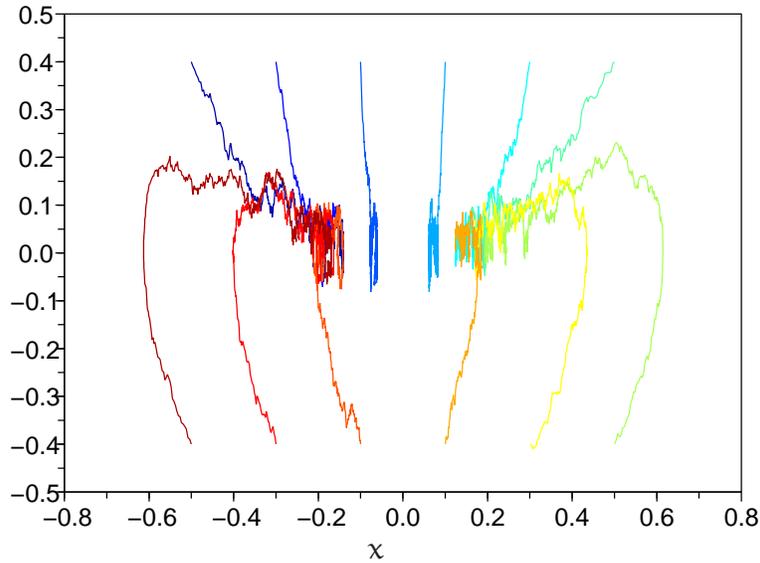} \\[-2ex]
        & $x$
    \end{tabular}
    \caption{sample trajectories of the example stochastic 
    system~\eqref{eq:toy} from different initial conditions
    for different realisations of the noise, $\sigma=0.05$\,.}
    \label{fig:sim1}
\end{figure}

Figure~\ref{fig:sim1} plots some typical trajectories of the \sde\ system~\eqref{eq:toy}.  In this domain near the origin the $y$~variable decays exponentially quickly to $y\approx x^2$\,; whereas the $x$~variable evolves relatively slowly over long times.  Thus the $y$~variable represents fast, microscopic, uninteresting modes, they are ``slaved'' to~$x$ \cite[]{Schoner86}, whereas the $x$~variable represents the long lasting, macroscopic modes of interest to the long term dynamics.  The white noise~$\sigma\phi(t)$ added to the $y$~dynamics induces fluctuations on all time scales.  We explore fundamental issues in the modelling of the multiscale dynamics of the \sde~\eqref{eq:toy}.

The challenge is to adapt the deterministic normal form transformation, Section~\ref{sec:ddd}, to the stochastic system~\eqref{eq:toy} in order to not only decouple the interesting slow modes, but to simplify them as far as possible, Section~\ref{sec:ssefp}. Two other alternative normal forms are explored in Sections \ref{sec:trsnfdm}~and~\ref{sec:aa}, but I argue that they are not so useful.  The analysis and argument is very detailed in order to demonstrate in a simple setting how Principles \ref{p:secular}--\ref{p:memory} are realised at the expense of having to anticipate future noise.  If you are familiar with the concept of stochastic normal forms, you could skip to Section~\ref{sec:nfsmsde} for generic arguments of the new results.

\subsection{Decouple the deterministic dynamics}
\label{sec:ddd}

\begin{figure}
    \centering
    \begin{tabular}{c@{}c}
        \raisebox{20ex}{$y$} & 
        \includegraphics{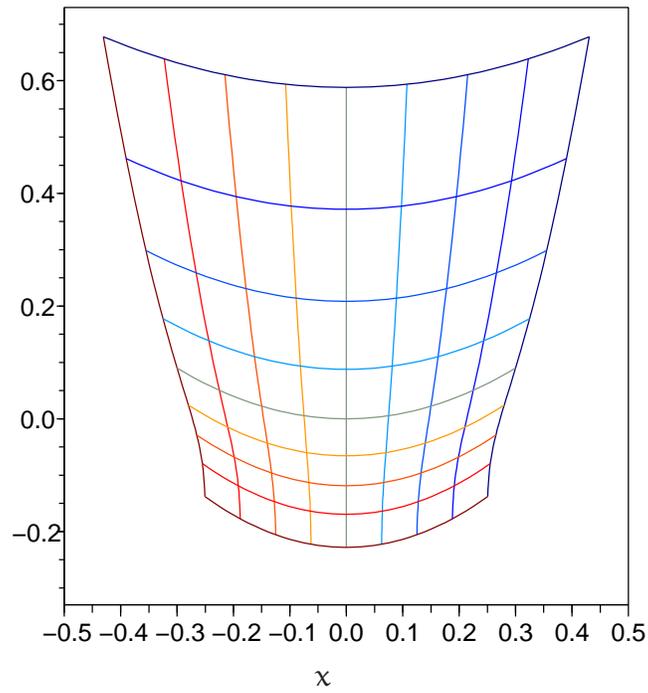} \\[-2ex]
        & $x$
    \end{tabular}
	\caption{coordinate curves in the $xy$-plane of the $(X,Y)$~coordinate
	system that simplifies to~\eqref{eq:toynf} the algebraic description
	of the dynamics of the deterministic ($\sigma=0$) system~\eqref{eq:toy}.}
    \label{fig:egnf}
\end{figure}

Initially consider the example toy system~\eqref{eq:toy} when there is no noise, $\sigma=0$\,.  A deterministic normal form coordinate transform decouples the deterministic slow and fast dynamics:
\begin{eqnarray}
    x&=&X +XY +\rat32XY^2 -2X^3Y +\rat52XY^3 +\cdots\,,
    \label{eq:toynfx}
    \\
    y&=&Y +X^2 +2Y^2 +4Y^3 -4X^2Y^2 +8Y^4 +\cdots\,.
    \label{eq:toynfy}
\end{eqnarray}
Figure~\ref{fig:egnf} shows the coordinate curves of this $(X,Y)$ coordinate system.  The coordinate transform is a near identity because near the origin $x\approx X$ and $y\approx Y$\,.  In the new $(X,Y)$  coordinate system, the evolution of the toy system~\eqref{eq:toy}
becomes
\begin{equation}
    \dot X= -X^3
    \quad\text{and}\quad
    \dot Y= -(1 +2X^2 +4X^4)Y +\cdots\,.
    \label{eq:toynf}
\end{equation}
Observe the $Y$-dynamics are that of exponentially quick decay to the slow manifold $Y=0$ at a rate $1+2X^2+4X^4+\cdots$\,. In the original variables, from the transform \eqref{eq:toynfx}~and~\eqref{eq:toynfy}, this slow manifold is the curve $x=X$ and $y=X^2$ \cite[]{Elphick87b}. The dynamics on this slow manifold, $\dot X=-X^3$ from~\eqref{eq:toynf}, form the accurate, macroscopic, long term model.

The slow $X$~dynamics are \emph{also} independent of the $Y$~variable and thus the initial value~$Y(0)$ and subsequent~$Y(t)$ are immaterial to the long term evolution. Thus to make accurate forecasts, project onto the slow manifold $Y=0$ along the coordinate curves of constant~$X$ \cite[]{Cox93b} seen in Figure~\ref{fig:egnf}. Equivalently, because the slow $X$~dynamics are independent of the $Y$~variable, the dynamics of the system~\eqref{eq:toy} map the curves of constant~$X$ in Figure~\ref{fig:egnf} into other curves of constant~$X$. Thus initial conditions on any one curve of constant~$X$ all evolve towards the same trajectory on the slow manifold.

But these comments are all for deterministic dynamics, $\sigma=0$\,. The next subsection answers the question: how can we adapt this beautifully simplifying coordinate transform to cater for stochastic dynamics?

\subsection{Simplify stochastic evolution as far as possible}
\label{sec:ssefp}

\begin{figure}
    \centering
    \begin{tabular}{c@{}c}
        \raisebox{25ex}{$y$} & 
        \includegraphics{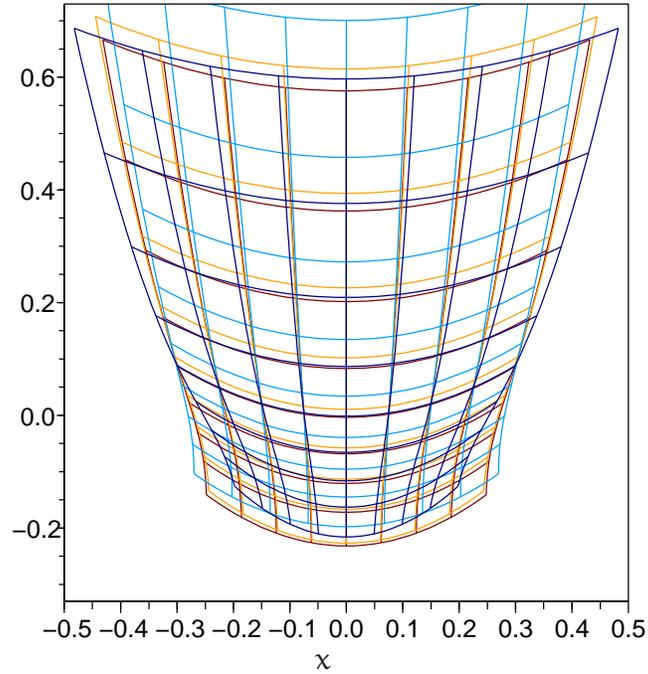} \\[-3ex]
        & $x$
    \end{tabular}
	\caption{four different (coloured) meshes represent \emph{either}
	four realisations sampled at one instant \emph{or} one realisation
	sampled at four instants of the coordinate curves in the
	$xy$-plane of the \emph{stochastic} $(X,Y)$~coordinate system that
	simplifies to~\eqref{eq:toy3}~and~\eqref{eq:toy4} the algebraic
	description of the dynamics of the stochastic ($\sigma=0.2$)
	system~\eqref{eq:toy}.}
    \label{fig:egsnf}
\end{figure}

Now explore the construction of a coordinate transform that decouples the fast and slow dynamics of the toy \sde~\eqref{eq:toy} in the presence of its stochastic forcing.  In order to cater for the stochastic fluctuations, the coordinate transform must be time dependent through dependence upon the realisation of the noise, as shown schematically in Figure~\ref{fig:egsnf}.  
This subsection is very detailed in order to argue that no alternatives go unrecognised.  The method is to iteratively refine the stochastic coordinate transform based upon the residuals of the governing toy \sde~\eqref{eq:toy}.  

Although our focus is
on the case when $\phi(t)$~is a white noise, because we use the usual
calculus of the Stratonovich interpretation, the algebraic results also
apply to smoother processes~$\phi(t)$.  For two examples, the
forcing~$\phi(t)$ could be the output of a deterministic chaotic
system~\cite[e.g.]{Just01}, or the forcing~$\phi(t)$ could be even as
regular as a periodic oscillator.  Thus the algebraic expressions
derived herein apply much more generally than to just white
noise~$\phi$.  However, the \emph{justification for the particular}
coordinate transform often depends upon the peculiar characteristics of
white noise.  For forcing~$\phi$ which is smoother than white noise,
although our results herein apply, other particular coordinate
transforms \emph{may be preferable} in order to achieve other desirable
outcomes in the transformation (outcomes not attainable when $\phi$~is
white noise).  These possibilities are not explored.  Instead, almost
everywhere throughout this article, the forcing~$\phi(t)$ denotes a
white noise process in a Stratonovich interpretation of \sde{}s.

Let us proceed to iteratively develop a stochastic coordinate transform of the \sde~\eqref{eq:toy} via stepwise refinement \cite[]{Roberts96a}.

\paragraph{First, consider the fast dynamics} With $x\approx X$ and $\dot X\approx 0$\,, seek a change to the $y$~coordinate of the form
\begin{equation}
    y= Y +\eta'(t,X,Y)+\cdots
    \quad\text{and}\quad
    \dot Y= -Y +G'(t,X,Y)+\cdots\,,
    \label{eq:yans1}
\end{equation}
where $\eta'$~and~$G'$ are small, $\Oeps2$, corrections to the transform and the corresponding evolution.  I introduce the parameter~$\epsilon$ to provide a convenient ordering of the terms that arise in the algebra: formally set $\epsilon=|(X,Y,\sigma)|$ with the effect that $\epsilon$~counts the number of~$X$, $Y$~and~$\sigma$ factors in any one term. Substitute~\eqref{eq:yans1} into the $y$ \sde~\eqref{eq:toy} and drop products  of small corrections to recognise we need to solve
\begin{equation}
    G'+\D t{\eta'}+\eta'-Y\D Y{\eta'}=\sigma\phi(t) +X^2-2Y^2 +\cdots\,;
    \label{eq:doty1}
\end{equation}
partial derivatives are here done keeping constant the other two variables of $X$, $Y$~and~$t$.

Consider first the deterministic terms.  To solve $G' +\eta'_t +\eta' -Y\eta'_Y=X^2-2Y^2$\,, keep the evolution as simple as possible (\P{ruthless}) by seeking corrections $G'=0$ and $\eta'=aX^2+bY^2$\,. Substitute into the equation to see $a=1$ and $b=2$\,.  Consequently include $\eta'=X^2+2Y^2$ into the coordinate transform.

Second, consider the stochastic term~$\sigma\phi(t)$ in the right-hand side of~\eqref{eq:doty1}.  To solve $G' +\eta'_t +\eta' -Y\eta'_Y =\sigma\phi(t)$ keeping the $Y$~dynamics as simple as possible (\P{ruthless}) choose the convolution~$\sigma\Z-\phi$\,, defined in~\eqref{eq:zmuf}, to be part of the correction~$\eta'$ to the coordinate transform.  Consequently the new approximation of the coordinate transform and the dynamics is
\begin{equation}
    y=Y+X^2+2Y^2+\sigma\Z-\phi+\cdots
    \quad\text{and}\quad
    \dot Y=-Y+\cdots\,.
\end{equation}
In these leading order terms in the coordinate transform,  see the stochastic slow manifold (\ssm) $Y=0$ corresponds to the vertically fluctuating parabola $y\approx X^2+\sigma\Z-\phi$ as seen in the overall vertical displacements of  the coordinate mesh in Figure~\ref{fig:egsnf}.

\paragraph{The convolution} For any non-zero parameter~$\mu$, and consistent with the convolution in the example transform~\eqref{eq:sde0nf}, define the convolution
\begin{equation}
    \Z{\mu}\phi=
    \begin{cases}
        \int_{-\infty}^t \exp[\mu(t-\tau)]\phi(\tau)\,d\tau\,,
        &\mu<0\,, \\
        \int_t^{+\infty} \exp[\mu(t-\tau)]\phi(\tau)\,d\tau\,,
        &\mu>0\,,             
    \end{cases}
    \label{eq:zmuf}
\end{equation}
so that the convolution is always with a bounded exponential (\P{secular}).  Five useful properties  of this convolution are
\begin{eqnarray}&&
    \Z\mu1=\frac1{|\mu|}\,,
    \label{eq:conv1}\\&&
    \frac{d\ }{dt}\Z{\mu}\phi=-\sgn\mu\,\phi+\mu\Z{\mu}\phi\,,
    \label{eq:ddtconv}
    \\&&
    E[\Z{\mu}\phi]=\Z{\mu}E[\phi]\,,
    \label{eq:exz}
    \\&&
    E[(\Z{\mu}\phi)^2]=\frac1{2|\mu|}\,,
    \label{eq:exzz}
    \\&&
    \Z\mu\Z\nu=\begin{cases}
    \frac1{|\mu-\nu|}\big[ \Z\mu+\Z\nu \big]\,, &\mu\nu<0\,, \\
    \frac{-\sgn\mu}{\mu-\nu}\big[ \Z\mu-\Z\nu \big]\,, 
    &\mu\nu>0\ \&\ \mu\neq\nu\,.
    \end{cases}
\end{eqnarray}
Also remember that although with $\mu<0$ the convolution~$\Z\mu$ integrates over the past, with $\mu>0$\,, as we will soon need, the convolution~$\Z\mu$ integrates into the future; both integrate over a time scale of order~$1/|\mu|$.

\paragraph{Second, consider the slow dynamics} Seek a correction to the
stochastic coordinate transform of the form
\begin{equation}
    x = X  +\xi'(t,X,Y)+\Oeps3
    \quad\text{and}\quad
    \dot X = F'(t,X,Y)+\Oeps3\,.
\end{equation}
where $\xi'$~and~$F'$ are~$\Oeps2$ corrections to the transform and the corresponding evolution.  Substitute into the $x$~equation of \sde~\eqref{eq:toy} and omit small products to recognise we need to solve
\begin{equation}
    F'+\D t{\xi'}-Y\D Y{\xi'}=-XY +\sigma X\Z-\phi+\Oeps3\,.
    \label{eq:dotx1}
\end{equation}
First, try $\xi'=aXY$ to find the deterministic term on the right-hand side is matched when
$a=1$\,.  
Second, consider the stochastic part of the equation:
$F'+{\xi'_t}-Y {\xi'_Y}=\sigma X\Z-\phi$\,.  The $\xi'_Y$ cannot help
us solve this stochastic part as there is no $Y$~factor in the
right-hand side term.  
We do not want to assign a fast time convolution
into the slow evolution~$F'$ (\P{memory}), but cannot integrate the
forcing~$\phi$ into~$\xi'$ as then terms would grow like the Wiener
process $W=\int\phi\,dt$ (\P{secular}).  
Instead, formally integrate by
parts to split $\Z-\phi=-\phi+\Z-\dot\phi$ and hence choose components
$F'=-\sigma X\phi$ and $\xi'=\sigma X\Z-\phi$\,.  Consequently
\begin{equation}
    x=X+XY+\sigma X\Z-\phi+\Oeps3
    \quad\text{and}\quad
    \dot X=-\sigma X\phi+\Oeps3\,.
    \label{eq:second}
\end{equation}

\paragraph{Third, reconsider the fast dynamics} Seek corrections,
$\eta'$~and~$G'$, to the $y$~transform and $Y$~evolution driven by
the updated residual of the $y$~equation in \sde~\eqref{eq:toy}:
\begin{equation}
    G'+\D t{\eta'} +\eta'-Y\D Y{\eta'} 
    = -4\sigma Y\Z-\phi  -2\sigma^2(\Z-\phi)^2 +\Oeps3\,.
\end{equation}
Separately consider the two stochastic forcing terms on the right-hand side.
\begin{itemize}
	\item To solve $G'+\eta'_t+\eta'-Y\eta'_Y = -4\sigma Y\Z-\phi$ we
	must seek $G'$~and~$\eta'$ proportional to~$Y$, whence
	$\eta'-Y\eta'_y=0$\,, and so to avoid secular terms in~$\eta'$
	(\P{secular}), and to avoid fast time convolution in the
	$Y$~evolution (\P{memory}), integration by parts enables to choose $G'=-4\sigma Y\phi$ and $\eta'=4\sigma Y\Z-\phi$\,.  
    
    \item To solve $G'+\eta'_t+\eta'-Y\eta'_Y = -2\sigma^2(\Z-\phi)^2$ seek
    $G'$~and~$\eta'$ independent of $X$~and~$Y$.  Hence choose $G'=0$
    (\P{ruthless}) and then the convolution
    $\eta'=-2\sigma^2\Z-(\Z-\phi)^2$ corrects the coordinate
    transform.\footnote{The right-hand side of this correction~$\eta'$ has
    non-zero mean.  We could assign the mean,~$-\sigma^2$, into the
    $Y$~evolution as a mean (downwards) forcing, but then this destroys
    $Y=0$ as being the slow manifold, contradicting \P{zero}.}

\end{itemize}

Consequently, the fast time transform and dynamics are more
accurately
\begin{eqnarray}
    y&=& Y +X^2 +2Y^2 +\sigma\big[\Z-\phi 
    +4 Y\Z-\phi \big] 
    \nonumber\\&&{}
    -2\sigma^2\Z-(\Z-\phi)^2 
    +\Oeps3\,,
    \\
    \dot Y&=&-Y -4\sigma Y\phi+\Oeps3\,.
\end{eqnarray}

\paragraph{Fourth, reconsider the slow dynamics}  Seek corrections to the
transform and evolution, $\xi'$~and~$F'$, driven by the
updated residual of the $x$~equation of \sde~\eqref{eq:toy}:
\begin{eqnarray}&&
    F'+\D t{\xi'}-Y\D Y{\xi'} = -X^3-3XY^2 +\sigma XY(5\phi-6\Z-\phi)
    \nonumber\\&&\quad{}
    +\sigma^2 X\left[ \phi\Z-\phi -(\Z-\phi)^2 +2\Z-(\Z-\phi)^2
    \right]
    +\Oeps4\,.
    \label{eq:res4x}
\end{eqnarray}
Consider the right-hand side term by term:
\begin{itemize}
    \item To account for the deterministic forcing, we must choose
    $F'=-X^3$ and $\xi'=\rat32XY^2$ in the traditional manner.
    
	\item To match the term linear in noise, $F'+\xi'_t-Y\xi'_Y = \sigma XY(5\phi-6\Z-\phi)$\,, we must seek $F'$~and~$\xi'$ proportional to~$XY$, whence
	$\xi'_t-Y\xi'_y\mapsto\xi'_t-\xi'$\,.  Consequently foreknowledge,
	anticipation, of the noise appears.  Consider the two cases:
    \begin{itemize}
		\item allowing anticipation (implementing \P{ruthless}) and in accord with \cite{Arnold98}, we
		assign all of this term to the coordinate transformation with
		$F'=0$ and $\xi'=-\sigma XY\Z+ (5\phi-6\Z-\phi)$\,;

		\item disallowing anticipation, we must assign all of this term
		into the $X$~evolution by assigning $F'= \sigma
		XY(5\phi-6\Z-\phi)$ and $\xi'=0$ ---the difficulty here being that
		the evolution to the \ssm\ then depends undesirably
		upon~$Y$, contradicting~\P{decouple}. Section~\ref{sec:aa}
        explores this case.
    \end{itemize}

	\item   For the quadratic noise term in~\eqref{eq:res4x}, seek contributions
	to the solution which are proportional to~$X$; consequently, on the
	left-hand side $-Y\xi_Y=0$\,.  Then we avoid secularity
	(\P{secular}) by extracting the mean of the right-hand side term
	and assign the mean into the evolution~$F'$; but at least part of the
	fluctuations cannot be assigned into the transform~$\xi'$ as the integral of noise is
	a Wiener process which almost surely is secular.
    
    Now, as in the earlier integration by parts, separate these quadratic noise terms into
    \begin{eqnarray*}
        (\Z-\phi)^2&=&\phi\Z-\phi-\rat12\rat{d}{dt}[ (\Z-\phi)^2]\,,
        \\
        \Z-(\Z-\phi)^2&=&(\Z-\phi)^2-\rat{d}{dt}[\Z-(\Z-\phi)^2]
        \\&=&\phi\Z-\phi-\rat{d}{dt}[ \rat12(\Z-\phi)^2
        +\Z-(\Z-\phi)^2]\,,
    \end{eqnarray*}
    and so these contribute corrections $F'=\sigma^2X\phi\Z-\phi$ and
    $\xi'=-\sigma^2 X(\rat12+2\Z-)(\Z-\phi)^2$\,.
\end{itemize}
The upshot is that the $x$~transformation and $X$~evolution is more
accurately
\begin{eqnarray}
    x&=&X +XY +\rat32XY^2 
    +\sigma\big[ X\Z-\phi +XY\Z+(-5\phi+6\Z-\phi) \big]
    \nonumber\\&&{}
    -\sigma^2X(\rat12+2\Z-)(\Z-\phi)^2
    +\Oeps4 \,,
    \label{eq:toyx}
    \\
    \dot X&=& -X^3 
    -\sigma X\phi
    +2\sigma^2X\phi\Z-\phi
    +\Oeps4\,.
    \label{eq:toy3}
\end{eqnarray}

\paragraph{Lastly, reconsider the fast dynamics}
Seek corrections $\eta'$~and~$G'$ to the $y$~transform and $Y$~evolution
driven by the updated residual of the $y$~equation in \sde~\eqref{eq:toy}:
\begin{eqnarray}&&
    G'+\D t{\eta'} +\eta'-Y\D Y{\eta'} 
    = -2X^2Y-8Y^3
    \nonumber\\&&\quad{}
    +\sigma\big[ 2 X^2(\phi-\Z-\phi)
    +8 Y^2(2\phi-3\Z-\phi)\big]
    \nonumber\\&&\quad{}
    +8\sigma^2 Y\big[ 2\phi\Z-\phi -2(\Z-\phi)^2 +\Z-(\Z-\phi)^2 \big]
    \nonumber\\&&\quad{}
    +8\sigma^3(\Z-\phi)\Z-(\Z-\phi)^2 \,.
\end{eqnarray}
Separately consider the forcing terms on the right-hand side.
\begin{itemize}
    \item The deterministic terms generate contributions $G'=-2X^2Y$ 
    and $\eta'=4Y^3$ in the usual manner.
    
    \item To solve $G'+\eta'_t+\eta'-Y\eta'_Y = +2\sigma
    X^2(\phi-\Z-\phi)$ we must seek $G'$~and~$\eta'$ proportional
    to~$X^2$ whence $-Y\eta'_y=0$\,.  Consequently, assign everything 
    to the transform employing the convolution~$\Z-$: $G'=0$ and $\eta'=2\sigma
    X^2(\Z-\phi-\Z-\Z-\phi)$ (\P{zero} and \P{ruthless}).

	\item The $Y^2$~term raises the issue of anticipation again.  To
	solve $G'+\eta'_t+\eta'-Y\eta'_Y =8\sigma Y^2(2\phi-3\Z-\phi)$ we
	must seek solutions proportional to~$Y^2$ and hence the left-hand
	side $G'+\eta'_t+\eta'-Y\eta'_Y=G'+\eta'_t-\eta'$\,.  By
	\P{ruthless}, assign the entire forcing into the transform by
	setting $G'=0$ and $\eta'=-8\sigma Y^2\Z+(2\phi-3\Z-\phi)$\,.  This
	requires anticipation of the forcing through $\Z+\phi$ and
	$\Z+\Z-\phi$\,.
    
    \item The remaining two terms do not generate new issues so I do 
    not describe the details.
    
\end{itemize}
The normal form coordinate transform for the  fast dynamics is thus
\begin{eqnarray}
    y&=& Y +X^2 +2Y^2 +4Y^3 
    +\sigma\big[\Z-\phi 
    +4 Y\Z-\phi 
    \nonumber\\&&\quad{}
    +2X^2\Z-(\phi-\Z-\phi)  
    -8Y^2\Z+(2\phi-3\Z-\phi) \big] 
    \nonumber\\&&{}
    +4\sigma^2 Y(1-2\Z-)(\Z-\phi)^2 
    \nonumber\\&&{}
    +8\sigma^3\Z-\big[(\Z-\phi)\Z-(\Z-\phi)^2 \big]
    +\Oeps4\,,
    \label{eq:toyy}
    \\
    \dot Y&=&-Y -2X^2Y -4\sigma Y\phi +8\sigma^2 Y\phi\Z-\phi +\Oeps4\,.
    \label{eq:toy4}
\end{eqnarray}

\paragraph{Higher order model}  Further algebra leads to the
construction of a stochastic coordinate transform from the original
$(x,y)$~variables to the new $(X,Y)$~variables so that the dynamics
of the example \sde~\eqref{eq:toy} is
\begin{eqnarray}
    \dot X&=& -X^3 
    -\sigma X\phi
    +2\sigma^2X\phi\Z-\phi
    \nonumber\\&&{}
    -4\sigma^2X^3\phi\Z-\Z-\phi
    +\Ord{\epsilon^6,\sigma^3}\,,
    \label{eq:toy6x}
    \\
    \dot Y&=&-(1 +2X^2 +4X^4)Y 
    -4\sigma(1+X^2) Y\phi 
    +8\sigma^2 Y\phi\Z-\phi 
    \nonumber\\&&{}
    +4\sigma^2 X^2Y\phi\big[ 3\Z-\phi -\Z+\phi -2\Z-\Z-\phi \big]
    \nonumber\\&&{}
    +\Ord{\epsilon^6,\sigma^3}\,.
    \label{eq:toy6y}
\end{eqnarray}
By employing a coordinate transform that depends upon the realisation of the noise, we maintain $Y=0$ as the exponentially attractive \ssm, see~\eqref{eq:toy6y},
independent of realisation.\footnote{Possibly, large enough noise~$\sigma$ could make $Y=0$ repulsive, just as large enough~$X$ may do.} 
However, we can only make the $X$~evolution independent of~$Y$, as
in~\eqref{eq:toy6x}, by anticipating noise, albeit anticipating only on
a fast time scale into the future.  
\emph{To rationally project onto the
\ssm\ we must accept some fast time scale anticipation.}

\paragraph{Irreducible fast time convolutions generate drift}
Also note that the $X$~and~$Y$ evolution equations,
\eqref{eq:toy6x}~and~\eqref{eq:toy6y}, contain algebraically
irreducible nonlinear noise such as $\phi\Z-\phi$, in defiance of
\P{memory}.  Over long times such irreducible noise could be replaced
by $\rat12+\rat1{\sqrt2}\alt\phi$ for some effectively new noise~$\alt\phi(t)$
\cite[]{Chao95}.  Such replacement was also justified by Khasminskii
(1996) as described by \cite{Srinamachchivaya90}.  Importantly, such
quadratic noise, in effect, generates a mean deterministic drift term
in the slow dynamics \cite[e.g.]{Srinamachchivaya90, Srinamachchivaya91}.
In applications such drifts can be vital.

\paragraph{The average SSM is not the deterministic slow manifold}
For the toy \sde~\eqref{eq:toy}, Section~\ref{sec:ddd} shows the deterministic
slow manifold is $y=x^2$\,.  In general the \ssm\ fluctuates about a
mean location which is different to this deterministic slow manifold.
From \eqref{eq:toyx}~and~\eqref{eq:toyy} with fast variable $Y=0$\,,
the \ssm\ is
\begin{eqnarray}
    x&=&X  
    +\sigma X\Z-\phi
    -\sigma^2X(\rat12+2\Z-)(\Z-\phi)^2
    +\Oeps4 \,,
    \label{eq:toyssm}
    \\
    y&=& X^2 
    +\sigma\big[\Z-\phi 
    +2X^2\Z-(\phi-\Z-\phi)  
    \big] 
    \nonumber\\&&{}
    +8\sigma^3\Z-\big[(\Z-\phi)\Z-(\Z-\phi)^2 \big]
    +\Oeps3\,,
\end{eqnarray}
Take expectations, and using \eqref{eq:exz}~and~\eqref{eq:exzz},
\begin{equation}
    E[x] =(1 -\rat52\sigma^2)X +\Ord{\epsilon^4}
    \qtq{and}
    E[y] =X^2 +\Ord{\sigma^3,\epsilon^4} \,.
    \label{eq:exy}
\end{equation}
Observe $E[y]\neq E[x]^2$\,, instead $E[y]\approx (1+\rat52\sigma^2)E[x]^2$ so that the average \ssm\ is a steeper parabola shape than the deterministic slow manifold.  
It is quadratic noise processes that deform the average \ssm\ from the deterministic.

\paragraph{Forecast from initial conditions}
Suppose at time $t=0$ we observe the state $(x_0,y_0)$, what forecast
can we make with the \ssm\ \sde~\eqref{eq:toy6x}?  Revert the
asymptotic expansion of the stochastic coordinate transform
\eqref{eq:toyx}~and~\eqref{eq:toyy} to deduce
\begin{eqnarray}
    X&=&x+x^3-xy+\rat32xy^2
    +2\sigma xy\Z+\phi
    -2\sigma^2 x(\Z+\phi)(\Z-\phi)
    \nonumber\\&&{}
    +\Ord{\epsilon^4}\,,
    \label{eq:toyrvx}\\
    Y&=&y-x^2-2y^2
    -\sigma\Z-\phi
    +\sigma^2(1+\Z-{})(\Z-\phi)^2
    +\Ord{\epsilon^3}\,.
    \label{eq:toyrvy}
\end{eqnarray}
Then the correct initial condition for the long term dynamics on the
\ssm, governed by the \sde~\eqref{eq:toy6x}, is the $X$~component of this
reversion,~\eqref{eq:toyrvx}, evaluated at the observed state, namely
\begin{equation}
    X(0)=x_0+x_0^3-x_0y_0+\rat32x_0y_0^2
    +2\sigma x_0y_0\Z+\phi
    -2\sigma^2 x_0(\Z+\phi)(\Z-\phi)
    +\Ord{\epsilon^4}\,.
    \label{eq:toyic}
\end{equation}
This is a projection of the observed initial state onto the \ssm\ to
provide an initial condition~$X(0)$ for the slow mode.  However, this
projection involves both memory and anticipatory convolutions of the
noise.  There are at least three interesting issues with computing this
initial~$X(0)$.  First, at the initial instant we do not know either
the future nor the past, so the terms involving the noise~$\phi$ are
unknown.  Using the expectations \eqref{eq:exz}~and~\eqref{eq:exzz},
the projection~$X(0)$ has known mean
\begin{displaymath}
    E[X(0)]=x_0+x_0^3-x_0y_0+\rat32x_0y_0^2
    +\Ord{\epsilon^4}\,,
\end{displaymath}
with known variance
\begin{displaymath}
    V[X(0)]\approx 2\sigma^2x_0^2y_0^2+\sigma^4x_0^2\,.
\end{displaymath}
That is, a given observed state $(x_0,y_0)$ corresponds to a stochastic
state for the evolution of the slow mode model on the \ssm.  Second,
but if this state~$X(0)$ for the slow mode is to be used in a
simulation to make forecasts of the future, then we know the future of
the noise~$\phi$.  The future values of noise~$\phi$ are just those we
use in integrating the slow mode \sde~\eqref{eq:toy6x}.  Thus for
simulation, we do eventually know the anticipatory convolutions~$\Z+\phi$
in~\eqref{eq:toyic}, but not the memory convolution~$\Z-\phi$.  In this
case the mean of the projection
\begin{displaymath}
    E[X(0)]=x_0+x_0^3-x_0y_0+\rat32x_0y_0^2+2\sigma x_0y_0\Z+\phi
    +\Ord{\epsilon^4}\,,
\end{displaymath}
with variance $V[X(0)]\approx 2\sigma^4x_0^2(\Z+\phi)^2$\,.  Lastly, if
we made additional observations for times $t<0$\,, then the additional
information could partially determine the past history of the
noise~$\phi$ and hence help us estimate the memory
convolution~$\Z-\phi$.  
These three cases emphasise that the initial state~$X(0)$ of the slow variable depends upon more than just the initial observed state~$(x_0,y_0)$.

\subsection{Try retaining some noise forcing of the decaying modes}
\label{sec:trsnfdm}

What if, contradicting both \P{ruthless} and \P{zero}, we allow some
forcing noise to remain in the $Y$~evolution?  Is there any useful
freedom?  Here I argue there is not.

\paragraph{First, consider the $y$~dynamics} Assume that we solve
equation~\eqref{eq:doty1} by choosing $\eta'=0$ and $G'=\sigma\phi$\,.
This keeps the noise imposed upon the rapid decay of the fast variable~$Y$.

\paragraph{Second, consider the $x$~dynamics}  The residual of the
$x$~equation in \sde~\eqref{eq:toy} then leads to solving
\begin{equation}
    F'+\xi'_t-Y\xi'_Y=-\sigma X\phi+\cdots\,,
    \label{eq:dotx2}
\end{equation}
instead of~\eqref{eq:dotx1}.  The solution must be proportional to~$X$,
so the $Y\xi'_Y$~term on the left-hand side is no use.  But secularity,
\P{secular}, implies we cannot assign any of the right-hand side component
into~$\xi$ and so we are forced to have $\xi'=0$ and $F'=-\sigma
X\phi$\,.  So far the coordinate transform itself remains identical to
the deterministic~\eqref{eq:toynfx}--\eqref{eq:toynfy}.

\paragraph{Third, reconsider the $y$~dynamics} Using the above
coordinate transform, determine more corrections from the new
residual by solving
\begin{equation}
    G'+\eta'_t+\eta'-Y\eta'_Y=-4\sigma Y\phi+2\sigma(X^2-6Y^2)\phi
    +\cdots\,.
\end{equation}
The terms on the right-hand side could all be assigned to the
$Y$~evolution as $G'=\sigma(-4Y +2X^2+4Y^2)\phi$ and $\eta'=0$\,.  The
term linear in~$Y$ certainly has to go into the $Y$~evolution, but the
others do not.  Instead, try an arbitrary convex combination of the above and
$G'=-4\sigma Y\phi +8\sigma^2Y\phi\Z+\phi$ (an anticipatory convolution
appears) and $\xi'=\sigma(2X^2\Z-\phi-4Y^2\Z+\phi)$\,.

\paragraph{Last, reconsider the $x$~dynamics} Seek corrections
$\xi'$~and~$F'$ to the $x$~transform and the $X$~evolution driven by the
new residual:
\begin{eqnarray}&&
    F'+\xi'_t-Y\xi'_Y = 2\sigma XY\phi +\cdots\,.
\end{eqnarray}
This residual is independent of the convex combination of the terms
just determined above by the $y$~dynamics.  The intractable difficulty
here is that the term~$2\sigma XY\phi$ on the right-hand side has to be assigned into~$F'$, contradicting \P{decouple}.  Further, we also get no hint of the
quadratic stochastic mean drift effect.  Since I parametrised the
freedom in the previous step, and the second step was forced, then to
obtain the quadratic mean drift term we must change the first step.
That is, we can only satisfy \P{decouple} and also can only extract the quadratic
mean drift term by abandoning the assumption that the first step is
useful.  Thus, allowing nonhomogeneous forcing of the $Y$~evolution is not useful for the \sde~\eqref{eq:toy}.

\subsection{Avoiding anticipation is less useful}
\label{sec:aa}

Alternatively, suppose we disallow anticipatory convolutions.  At least
one of the Principles~\ref{p:secular}--\ref{p:memory} then has to be
abandoned.  Principles \ref{p:ruthless}~and~\ref{p:memory} are ``as
possible'' principles so we meet them as best we can.  Well ordered
asymptotic expansions are essential, so avoid secularity, \P{secular}.
Section~\ref{sec:trsnfdm} shows abandoning \P{zero} is ineffective.
Consequently, in this subsection we explore abandoning \P{decouple},
the requirement to completely decouple the slow modes from the fast
modes.  But abandoning this principle means we are no longer able to use
the slow model to make high accuracy forecasts from every initial
condition.

The construction of a normal form transform that avoids anticipatory
convolutions follows the same sort of steps as described in detail in
Section~\ref{sec:ssefp}.  There is no point in redoing such detail.
Instead I state the coordinate transform as derived and checked by
computer algebra that also derived and checked the transform leading to
\eqref{eq:toy6x}~and~\eqref{eq:toy6y}.

The stochastic coordinate transform
\begin{eqnarray}
    x&=&X +\sigma X \Z-\phi  -\sigma^2 X(\rat12+2\Z-{})(\Z-\phi)^2 
    +\Ord{\epsilon^4}\,,
    \label{eq:aax}
    \\
    y&=&Y +X^2 +\sigma \big[1+4Y +2X^2(1-\Z-{})\big]\Z-\phi
    \nonumber\\&&{}
    +\sigma^2 \big[-2\Z-{} +4Y(1-2\Z-{}) \big](\Z-\phi)^2
    +\Ord{\epsilon^4}\,,
    \label{eq:aay}
\end{eqnarray}
transforms the system of \sde{}s~\eqref{eq:toy} to the equivalent system
\begin{eqnarray}
    \dot X&=& -X^3 -XY -\sigma X\phi -4\sigma XY\Z-\phi
    +2\sigma^2X\phi\Z-\phi +\Ord{\epsilon^4}\,,
    \label{eq:aaxx}
    \\
    \dot Y&=& -(1+2X^2+2Y)Y -4\sigma Y(\phi+2Y\Z-\phi)
    +8\sigma^2Y\phi\Z-\phi
    \nonumber\\&&{}
    +\Ord{\epsilon^4}\,.
    \label{eq:aayy}
\end{eqnarray}
Note two aspects: there are no anticipatory convolutions; and $X$~and~$Y$ variables are different to those of Section~\ref{sec:ssefp}.

\paragraph{The stochastic slow manifold is attractive}
This \sde\ system has $Y=0$ as an invariant manifold as every term in the
$Y$ \sde~\eqref{eq:aayy} is multiplied by~$Y$.  Thus, in the
transformed coordinates the stochastic slow manifold (\ssm) is $Y=0$\,.
This \ssm\ exponentially quickly attracts at least some finite domain
about the origin in $(X,Y,\sigma)$~space as the dominant terms in the
evolution are $\dot Y\approx -Y$\,.  From the stochastic coordinate
transform~\eqref{eq:aay}, the \ssm\ is parametrically given by~\eqref{eq:aax}
and
\begin{displaymath}
    y=X^2 +\sigma \big[1 +2X^2(1-\Z-{})\big]\Z-\phi
    -2\sigma^2 \Z-\big[(\Z-\phi)^2\big]
    +\Ord{\epsilon^4}\,.
\end{displaymath}
On this \ssm\ the evolution, from~\eqref{eq:aaxx}, is identical to the
\ssm\ model~\eqref{eq:toy6x}.  This normal form coordinate transform easily
displays the \ssm.

\paragraph{We cannot make accurate forecasts}
Suppose we specify some initial state~$(X_0,Y_0)$, either deterministic
or stochastic.  What forecast can we easily make with the \ssm\
model~\eqref{eq:toy6x}?  In general, none.  The reason is that in the
evolution to the \ssm, the \sde~\eqref{eq:aaxx} shows the slow
$X$~dynamics are coupled to the fast $Y$~dynamics.  But the point of
deriving a slow model, for most purposes, is to avoid resolving the
details of the fast dynamics; thus we cannot rationally project
from~$(X_0,Y_0)$ onto the \ssm.  In contrast, the normal form of a
deterministic system empowers rational projection from nearby initial
conditions onto the slow model for accurate forecasts \cite[]{Cox93b}.
Abandoning \P{decouple} means we cannot make accurate forecasts.

Because it adheres to \P{decouple}, the stochastic normal form of
Section~\ref{sec:ssefp}, similarly to the deterministic normal form,
empowers rational projection from nearby initial conditions onto the
\ssm.  But there is a catch: in order to do the projection we
need to anticipate the noise.  Since we generally will not know the
future noise, the stochastic normal form of Section~\ref{sec:ssefp}
also cannot be used for accurate forecasting.  In this sense the two
stochastic normal forms have equivalent power.  However, there is a
difference.  The anticipatory stochastic normal form of
Section~\ref{sec:ssefp} has explicit convolutions for the projection:
we may not know what they are, but we could certainly use the
convolutions to estimate bounds and distributions for the projection of
initial conditions.  In contrast, the stochastic normal form of this
section keeps the such information encrypted in the coupled fast and
slow dynamics of \eqref{eq:aaxx}~and~\eqref{eq:aayy}.  Consequently,
\emph{maintaining \P{decouple}, decoupling the slow modes from the
fast, appears more powerful than avoiding anticipatory convolutions.}

%
%

\section{Normal forms of SDEs for long term modelling}
\label{sec:nfsmsde}

This section uses formal arguments to establish a couple of key generic
properties of stochastic normal forms seen in the example \sde\ system
of the previous section.  We establish firstly that a stochastic
coordinate transform can decouple slow modes from fast modes, to make
the stochastic slow manifold (\ssm) easy to see, and secondly that
although anticipation of the noise may be necessary in the full
transform no anticipation need appear on the \ssm.

Consider a general system of \sde{}s for variables $\vec
x(t)\in\mathbb R^m$ and $\vec y(t)\in\mathbb R^n$\,:
\begin{eqnarray}
    \dot{\vec x}&=&A\vec x+\vec f(\vec x,\vec y,t)\,,
    \label{eq:sdex}\\
    \dot{\vec y}&=&B\vec y+\vec g(\vec x,\vec y,t)\,,
    \label{eq:sdey}
\end{eqnarray}
where 
\begin{itemize}
	\item the spectrum of~$A$ is zero and for simplicity we assume
	$A$~is upper triangular with all elements zero except possibly
	$A_{i,j}$ for $j>i$ (such as in the Jordan form appropriate for
	position and velocity variables of a mechanical system); 
    
    \item for simplicity assume matrix~$B$ has been diagonalised with
    diagonal elements $\beta_1,\ldots,\beta_n$, possibly complex, with
    $\Re\beta_j<0$\,;\footnote{If matrix~$B$ is in Jordan form, rather than
    diagonalisable, then extensions of the arguments lead to the same
    results.}

    \item  $\vec f$ and $\vec g$ are stochastic functions that are
    ``nonlinear'', that is, $f$~and~$g$ and their gradients in $\vec 
    x$~and~$\vec y$ are all zero at the origin;
    
	\item the stochastic nature of the system of \sde{}s arises through
	the dependence upon the time~$t$ in the nonlinearity $\vec
	f$~and~$\vec g$---assume the time dependence is implicitly due to
	some number of independent white noise processes~$\phi_k(t)$ (which
	are derivatives of independent Wiener processes).
    
\end{itemize}
For such systems, \cite{Boxler89} guarantees the existence, relevance
and approximability of a stochastic centre manifold
for~(\ref{eq:sdex}--\ref{eq:sdey}) in some finite neighbourhood of the
origin.  We call this a stochastic slow manifold (\ssm) because we
assume matrix~$A$ does not have complex eigenvalues (oscillatory dynamics are considered briefly in Section~\ref{sec:osc}).

For example, the toy \sde\ system~\eqref{eq:toy} takes
the form~(\ref{eq:sdex}--\ref{eq:sdey}) with variables $\vec
x=(\sqrt\sigma,x)$ and $\vec y=y$\,, then
\begin{displaymath}
    A=
    \begin{bmatrix}
        0&0\\0&0
    \end{bmatrix} ,\quad
    B=-1\,,\quad
    \vec f=
    \begin{bmatrix}
        0\\-xy
    \end{bmatrix} ,\quad
    \vec g=x^2-2y^2+\sqrt\sigma^2\phi(t)\,.
\end{displaymath}

In principle, the matrices $A$~and~$B$ could also depend upon
the realisation of the noise.  When the Lyapunov exponents of the
corresponding linear dynamics are zero and negative respectively,
then a stochastic centre manifold still exists and has nice
properties \cite[]{Boxler89}.  However, here I restrict attention to
the algebraically more tractable case when the basic linear operators
$A$~and~$B$ are deterministic.

\paragraph{Stochastic singular perturbation systems} such as those
explored by \cite{Berglund03}, are a subset of the systems encompassed
by~(\ref{eq:sdex}--\ref{eq:sdey}).  For example,  transform the
deterministic singular perturbation system
\begin{equation}
    \dot x=f(x,y)\,,\quad \dot y=\frac1\epsilon g(x,y)\,,
\end{equation}
into the form~(\ref{eq:sdex}--\ref{eq:sdey}).  First, change to the
fast time $\tau=t/\epsilon$ so that
\begin{displaymath}
    \frac{dx}{d\tau}=\epsilon f(x,y)\,,\quad
    \frac{dy}{d\tau}=g(x,y)\,.
\end{displaymath}
Then change to a coordinate system $\xi$~and~$\eta$, where $\eta=0$
is the curve $g(x,y)=0$\,, in which the system
takes the form
\begin{displaymath}
    \frac{d\xi}{d\tau}=\epsilon F(\xi,\eta)\,,\quad
    \frac{d\eta}{d\tau}=\beta(\xi)\eta+G(\xi,\eta)\,.
\end{displaymath}
Consequently, in variables $\vec x=(\sqrt\epsilon,\xi)$ and $\vec
y=\eta$\,, the curve $(\vec x,\vec y)=(0,\xi,0)$ are a set of
equilibria, at each of which the dynamics are of the
form~(\ref{eq:sdex}--\ref{eq:sdey}).  Consequently there exists a slow
manifold around each point of the curve, which as a whole forms a slow
manifold in a neighbourhood of the curve \cite[]{Carr81}.  Thus, the
analysis presented here also applies to singular perturbation problems
by a change in time scale and coordinate system.

\paragraph{A stochastic coordinate transform}

We  transform the \sde~(\ref{eq:sdex}--\ref{eq:sdey}) in
$(\vec x,\vec y)$ to the new $(\vec X,\vec Y)$ coordinate system
by a stochastic, near identity, coordinate transform
\begin{equation}
    \vec x=\vec X+\vec \xi(\vec X,\vec Y,t)
    \quad\text{and}\quad
    \vec y=\vec Y+\vec\eta(\vec X,\vec Y,t)\,.
    \label{eq:xform}
\end{equation}
This stochastic coordinate transform is to be chosen such that the
\sde~(\ref{eq:sdex}--\ref{eq:sdey}) transforms to a ``simpler'' form
from which we may easily extract the \ssm.  Based upon the experience
of Section~\ref{sec:ssefp}, we seek to simplify the \sde{}s according
to Principles~{p:secular}--\ref{p:memory}, and allowing anticipation.

\subsection{Transform the fast dynamics}

Suppose~\eqref{eq:xform} is some approximation to the desired
coordinate transform.  Iteratively we seek corrections
$\vec\xi'$~and~$\vec\eta'$ to the transform, namely
\begin{equation}
    \vec x=\vec X+\vec \xi(\vec X,\vec Y,t)+\vec\xi'(\vec X,\vec Y,t)
    \quad\text{and}\quad
    \vec y=\vec Y+\vec\eta(\vec X,\vec Y,t)+\vec\eta'(\vec X,\vec Y,t)\,.
    \label{eq:xformd}
\end{equation}
Find corrections such that the corresponding updates to the
evolution, say $\vec F'$~and~$\vec G'$ in
\begin{eqnarray}
    \dot{\vec X}&=&A\vec X+\vec F(\vec X,\vec Y,t)
    +\vec F'(\vec X,\vec Y,t)\,,
    \label{eq:sdexxd}\\
    \dot{\vec Y}&=&B\vec Y+\vec G(\vec X,\vec Y,t)
    +\vec G'(\vec X,\vec Y,t)
    \label{eq:sdeyyd}\,,
\end{eqnarray}
are as simple as possible (\P{ruthless}).  

For the fast dynamics, the iteration is to substitute the corrected
transform~\eqref{eq:xformd} and
evolution~\eqref{eq:sdexxd}--\eqref{eq:sdeyyd} into the governing
\sde~\eqref{eq:sdey} for the fast variables.  Then drop products of
corrections as being negligible, and approximate coefficients of
corrections by their leading order term.  Then the equation for the
$j$th~component of the correction to the transform of the fast variable
and the new fast dynamics is
\begin{equation}
    G'_j+\D t{\eta'_j}-\beta_j\eta'_j +\sum_{\ell=1}^n \beta_\ell
    Y_\ell\D{Y_\ell}{\eta'_j} =\res_{\protect\ref{eq:sdey},j}\,.
    \label{eq:yd}
\end{equation}
Here, $\res_{\protect\ref{eq:sdey},j}$ denotes the residual of the
$j$th~component of the \sde~\eqref{eq:sdey}.  In constructing a
coordinate transform we repeatedly solve equations of this form to find
corrections.

We find what sort of terms may be put into the transformation~$\eta$
and what terms have to remain in the $\vec Y$~evolution by considering
the possibilities for the right-hand side.  The transform is
constructed as a multivariate asymptotic expansion about the origin in $(\vec X,\vec Y)$~space; thus all terms
in the \sde{}s are correspondingly written as asymptotic expansions
in~$(\vec X, \vec Y)$.  Suppose the right-hand side, the
residual~$\res_{\protect\ref{eq:sdey},j}$, has, among many others, a term of
the multinomial form
\begin{displaymath}
    c(t)\vec X^{\vec p}\vec Y^{\vec q}
    =c(t)\prod_{i=1}^mX_i^{p_i}
    \prod_{j=1}^nY_j^{q_j}\,,
\end{displaymath}
for some vectors of integer exponents $\vec p$~and~$\vec q$.  Because
of the special form of the `homological' operator on the left-hand
side of~\eqref{eq:yd}, seek contributions to the corrections of
\begin{displaymath}
    G'_j=a(t)\vec X^{\vec p}\vec Y^{\vec q}
    \quad\text{and}\quad
    \eta'_j=b(t)\vec X^{\vec p}\vec Y^{\vec q}\,.
\end{displaymath}
Then this component of~\eqref{eq:yd} becomes
\begin{equation}
    a+\dot b -\mu b=c
    \quad\text{where}\quad
    \mu= \beta_j-\sum_{\ell=1}^n q_\ell\beta_\ell\,.
    \label{eq:mud}
\end{equation}
Three cases arise.
\begin{enumerate}
	\item In the case $\mu=0$\,, we need to solve $a+\dot b=c$ where we
	want to put as much into~$b$ as possible (\P{ruthless}).
	Generically, the forcing~$c(t)$ will have mean and stochastically
	fluctuating components.  Neither of these can be integrated
	into~$b$ as they both give rise to secular terms (\P{secular}): the
	mean of~$c$ generates linear growth; the stochastically fluctuating
	part of~$c$ almost surely generates square-root
	growth.\footnote{\label{fn:det}In contrast, when the
	forcing~$c(t)$ is periodic, instead of stochastic, then the the
	forcing~$c(t)$ may be integrated into the coordinate transform~$b$,
	instead of being assigned to the evolution~$a$.} Thus the generic
	solution is $a=c$ and $b=0$\,, that is, assign $c(t)\vec X^{\vec p}
	\vec Y^{\vec q}$ to the $\vec Y$~evolution and nothing into the
	coordinate transform~$\vec\eta$.
    
	One example when this case occurs is when all exponents
	in~$\vec q$ are zero except for $q_j=1$\,.  Then the contribution
	to the $\vec Y$~evolution is simply $c(t)\vec X^{\vec p}Y_j$.
	Being linear in~$Y_j$, this contribution maintains \P{zero} that
	$\vec Y=\vec 0$ is the \ssm\ in the transformed coordinates.  
    
	In general, since~$\Re\beta_\ell$ are \emph{all}
	negative,\footnote{More generally, provided the
	eigenvalues~$\beta_j$ are all non-zero whether real or complex, the
	argument still holds.  Thus we can maintain \P{zero} over a very wide
	range of circumstances.} the case $\mu=0$ can only arise when at
	least one of the exponents~$\vec q$ is positive in order for the
	sum in~\eqref{eq:mud} to be zero.  Hence, generally there will be
	at least one $Y_\ell$~factor in updates~$\vec G'$ to the $\vec
	Y$~evolution, and so we maintain that $\vec Y=\vec 0$ is the \ssm.

	\item When $\Re\mu<0$\,, a solution of~\eqref{eq:mud} is to place all
	the forcing into the \ssm, $b=\Z\mu c$\,, and do not introduce a
	component into the $\vec Y$~evolution, $a=0$\,.  As $\Re\mu<0$\,, the
	convolution is over the past history of the noise affected
	forcing~$c(t)$; the convolution represents a memory of the forcing
	over a time scale of~$1/|\Re\mu|$.  This case of $\Re\mu<0$ arises when
	$-\Re\beta_j$ is large and the exponents~$\vec q$ are relatively
	small, corresponding to low order nonlinear factors of a rapidly
	dissipating mode.
    
	However, for large enough exponents~$\vec q$, that is for high
	enough order nonlinear terms, the rate~$\Re\mu$ must eventually become
	positive.  In the transition from negative to positive, the
	rate~$\Re\mu$ may become close to zero.  Then the time scale~$1/|\Re\mu|$
	becomes large and may be as large as the macroscopic time scale of
	the slow dynamics of interest.  In that case set the transform $b=0$ and assign
	this term in the forcing into the $\vec Y$~evolution with $a=c$\,.  The intended use of a macroscopic model defines a slow time scale and consequently affects which terms appear in the model.

	\item When $\Re\mu>0$\,, and accepting anticipation in the transform,
	we simply set $b=\Z\mu c$\,, and do not change the $\vec
	Y$~evolution, $a=0$\,.
    
\end{enumerate}
Consequently, \emph{we are always able to find a coordinate transform
which maintains that $\vec Y=\vec 0$ is the \ssm.}

You may have noticed that I omitted a term in~\eqref{eq:yd}: the term
$\D{X_\ell}{\eta'_j}A_{\ell,i}X_i$ should perhaps appear in the
left-hand side.  However, my omission is acceptable when the matrix~$A$
is upper triangular as then any term introduced which involves~$X_\ell$
only generates extra terms which are lower order in~$X_\ell$.  Although
such extra terms increase the order of~$X_i$ for $i>\ell$\,,
successive iterations generate new terms involving only fewer factors
of~$X_\ell$ and so iteration steadily accounts for the introduced
terms.  Similarly for the $\vec Y$~variables when the linear
operator~$B$ is in Jordan form due to repeated eigenvalues.  Discussing
equation~\eqref{eq:yd} for corrections is sufficient.  Analogous
comments apply to the the slow dynamics to which we now turn.

\subsection{Transform the slow dynamics}

For the slow dynamics, each iteration towards constructing a stochastic
coordinate transform substitutes corrections to the
transform~\eqref{eq:xformd} and the
evolution~(\ref{eq:sdexxd}--\ref{eq:sdeyyd}) into the governing
\sde~\eqref{eq:sdex} for the slow variables.  Then drop products of
corrections as being negligible, and approximate coefficients of
corrections by their leading order term.  Then the equation for the
$j$th~component of the correction to the transform of the slow variable is
\begin{equation}
    F'_j+\D t{\xi'_j} +\sum_{\ell=1}^n \beta_\ell
    Y_\ell\D{Y_\ell}{\xi'_j} =\res_{\protect\ref{eq:sdex},j}\,.
    \label{eq:xd}
\end{equation}
Here, $\res_{\protect\ref{eq:sdex},j}$ denotes the residual of the
$j$th~component of the \sde~\eqref{eq:sdex} evaluated at the current
approximation.  In constructing a stochastic coordinate transform we
repeatedly solve equations of this form to find corrections.  The
difference with the previous discussion of the fast variables is that
the left-hand side of~\eqref{eq:xd} does not have an analogue of the
$-\beta_j\eta'_j$~term.

We find what sort of terms may be put into the transformation
correction~$\vec \xi'$ and what terms have to remain in the $\vec
X$~evolution, via the correction~$\vec F'$, by considering the range of
possibilities for the right-hand side.  In general, the right-hand side
residual~$\res_{\protect\ref{eq:sdex},j}$ is a sum of terms of the form
\begin{displaymath}
    c(t)\vec X^{\vec p}\vec Y^{\vec q}
    =c(t)\prod_{i=1}^mX_i^{p_i}
    \prod_{j=1}^nY_j^{q_j}\,,
\end{displaymath}
for some vectors of integer exponents $\vec p$~and~$\vec q$.  Because
of the special form of the `homological' operator on the left-hand
side of~\eqref{eq:xd}, seek corresponding corrections
\begin{displaymath}
    F'_j=a(t)\vec X^{\vec p}\vec Y^{\vec q}
    \quad\text{and}\quad
    \xi'_j=b(t)\vec X^{\vec p}\vec Y^{\vec q}\,.
\end{displaymath}
Then~\eqref{eq:xd} becomes
\begin{equation}
    a+\dot b -\mu b=c
    \quad\text{where}\quad
    \mu= -\sum_{\ell=1}^n q_\ell\beta_\ell\,.
    \label{eq:mux}
\end{equation}
Two  cases typically arise.\footnote{The case $\Re\mu<0$ cannot arise as all the decay rates~$-\Re\beta_j>0$ when there are no fast unstable modes.  }
\begin{enumerate}
	\item The case $\mu=0$ only arises when the exponents $\vec q=\vec
	0$ as the exponents have to be non-negative and $\Re\beta_\ell<0$\,.  We need to solve
	$a+\dot b=c$ where we want to put as much into~$b$ as possible
	(\P{ruthless}).  Generically, the forcing~$c(t)$ will have mean and
	stochastically fluctuating components.  Neither of these can be
	integrated into~$b$ as they both give rise to secular terms
	(\P{secular}): the mean of~$c$ generates linear growth; the
	fluctuating part of~$c$ almost surely generates square-root
	growth.\footnote{But see footnote~\ref{fn:det},
	p.\pageref{fn:det}.} Thus at first sight the generic solution is
	$a=c$ and $b=0$\,, that is, assign $c(t)\prod_{i=1}^mX_i^{p_i}$ to
	the $\vec X$~evolution and nothing into the coordinate transform.
    
    But recall \P{memory}: we do not want fast time integrals in
    the slow evolution.  Consider the case when the forcing has the
    form of a fast time convolution $c=\Z\nu C$\,. 
    From~\eqref{eq:ddtconv} deduce
    \begin{displaymath}
        c=\Z\nu C=\frac1{|\nu|}C+\frac1\nu\frac{d\ }{dt}(\Z\nu C)
        =\frac1{|\nu|}C+\frac1\nu\frac{dc}{dt}\,.
    \end{displaymath}
	Hence to avoid fast time memory integrals in the slow $\vec
	X$~evolution (\P{memory}), set $a=C/{|\nu|}$ and $b=c/\nu$\,.
	If $C(t)$~in turn is a fast time convolution, then continue the above
	separation.  This separation corresponds to the integration by
    parts that Section~\ref{sec:ssefp} uses to avoid fast time,
    memory convolutions in the slow evolution.
    
	When the forcing~$c$ is a quadratic product of convolutions, then
	similar transformations and integration by parts eliminates all
	memory from the slow variables except terms of the form $\alt c(t)\Z\nu
	\aalt c(t)$ where $\alt c$~has no convolutions.  Algebraic transformations
	cannot eliminate such terms; for now accept the violation of
	\P{memory} in such quadratic forcing terms.

	Since the case $\mu=0$ can only arises for terms in the residual
	with no $\vec Y$~dependence, we can maintain that the slow
	evolution of the $\vec X$~variables are independent of~$\vec Y$,
	and this holds both on and off the \ssm.\footnote{However, when the
	fast dynamics contain rapidly oscillating, non-decaying modes, the
	corresponding eigenvalues occur as complex conjugate pairs which
	typically interact to cause $\mu=0$\,; among rapid oscillation we
	cannot completely decouple the slow modes from the fast oscillations
	\cite[]{Cox93b}.  Physically, waves do interact with mean flow.}
    
	\item The remaining case when $\Re\mu>0$ occurs when at least one of the exponents~$\vec q$ is positive.
	Accepting anticipation in the transform, we simply assign $b=\Z\mu
	c$\,, and do not change the $\vec X$~evolution, $a=0$\,.
\end{enumerate}
By anticipating noise we are \emph{always able to find a coordinate
transform which maintains a slow $\vec X$~evolution that is independent of
whether the system is on or off the \ssm.} Thus the projection of
initial conditions and the exponential approach to a solution of the
slow variables, called asymptotic completeness by \cite{Robinson96}, is
only assured by anticipation of the noise.

The preceding arguments are phrased in the context of an iteration
scheme to construct the stochastic coordinate transform and the
corresponding evolution.  Each step in the iterative process satisfies
the governing \sde{}s to higher order in the asymptotic expansions.  By
induction, we immediately deduce the following proposition.

\begin{proposition}
    \label{thm:nf}
with stochastic anticipation allowed, a near identity stochastic
coordinate transformation exists to convert the stochastic
system~(\ref{eq:sdex}--\ref{eq:sdey}) into the normal form
\begin{eqnarray}
    \dot{\vec X}&\simeq&A\vec X+\vec F(\vec X,t)\,,
    \label{eq:sdenfxx}\\
    \dot{\vec Y}&\simeq&\big[ B+G(\vec X,\vec Y,t) \big]\vec Y\,,
    \label{eq:sdenfyy}
\end{eqnarray}
where $\simeq$~denotes that these are equalities to any power of $(\vec
X,\vec Y)$ in an asymptotic expansion about the origin; there generally
are exponentially small errors.

Note: $\vec F$~and~$G$ may contain fast
time memory integrals but these need only occur as products with other
noise processes; for example, see \eqref{eq:toy3}~and~\eqref{eq:toy4}.
\end{proposition}

This proposition corresponds to the general Theorem~2.1 of
\cite{Arnold93} except they do not identify that memory
integrals may be mostly eliminated.

\subsection{Slow dynamics do not need to anticipate the noise}
\label{sec:sddnntatn}

Despite the presence of anticipatory convolutions appearing in the
stochastic coordinate transform, we here argue that none of them appear
in the slow dynamics because the anticipatory convolutions always involve
fast variables.  \cite{Bensoussan95} correspondingly show we do not need to anticipate
noise on a stochastic inertial manifold.

In the previous sections, the anticipatory convolutions only occur when
the rate~$\mu>0$\,.  But for both the slow and the fast components,
this occurrence can only be generated when at least one fast
$Y_j$~variable appears in the term under consideration.  Moreover,
there is no ordinary algebraic operation that reduces the number of
$\vec Y$~factors in any term: potentially the time derivative operator
might,
\begin{displaymath}
    \frac{d\ }{dt}=\D t{}+\sum_{\ell,k}Y_kB_{\ell,k}\D {Y_\ell}{}
    +\sum_{\ell,k}X_kA_{\ell,k}\D {X_\ell}{}\,,
\end{displaymath}
but although for non-diagonal $A$~and~$B$,  in the algebra $X_\ell$~variables may be
replaced by~$X_k$ and $Y_\ell$~variables may be replaced by~$Y_k$,
nonetheless the same number of variables are retained in each term and
a $\vec Y$~variable is never replaced by an $\vec X$~variable.  The
reason is that the $x$~and~$y$ dynamics are \emph{linearly} decoupled
in the system~(\ref{eq:sdex}--\ref{eq:sdey}).  Consequently all
anticipatory convolutions appear in terms with at least one component
of the fast variables~$\vec Y$.

Since the evolution~\eqref{eq:sdenfxx} of the slow modes~$\vec X$ is
free of $\vec Y$~variables, the evolution is also free of anticipatory
convolutions.  However, as seen in examples, there generally are
anticipatory convolutions in the evolution~\eqref{eq:sdenfxx} of the
fast modes~$\vec Y$.  Further, although the stochastic coordinate
transform~\eqref{eq:xform} has anticipatory convolutions, on the \ssm\ $\vec Y=\vec 0$ there are none.  Consequently the preceding
formal analysis leads to the following proposition.

\begin{proposition}
    \label{thm:mem}
although stochastic anticipation may be invoked, there \emph{need not}
be any anticipation in the dynamics~\eqref{eq:sdenfxx} of the slow
modes in the stochastic normal form of the
system~(\ref{eq:sdex}--\ref{eq:sdey}).
Moreover, on the \ssm, $\vec Y\simeq\vec 0$\,, the stochastic
coordinate transform~\eqref{eq:xform} need not have anticipation.
\end{proposition}

In contrast, \cite{Arnold93} and \cite{Arnold98} record anticipatory
convolutions in the slow modes of their examples, respectively (12)~and~(4.6).  Such
anticipatory convolutions are undesirable in using the normal form to
support macroscale models.

\section{Implications for multiscale modelling}
\label{sec:imstm}
This section describes some of the generic consequences of the previous
sections in modelling stochastic systems.

\paragraph{Anticipation}
All who write down and then use coarse scale models of
stochastic dynamics implicitly are soothsayers.
In writing down a coarse scale model, researchers neglect
the many details of any quickly decaying insignificant ignored modes.
Proposition~\ref{thm:nf} assures us that normally this neglect requires
us to know aspects of the near future of the ignored modes in order to
decouple the coarse modes from the uninteresting details.
In particular, providing initial conditions for the coarse model
requires looking into the future.
Nonetheless, Proposition~\ref{thm:mem} assures us that non-anticipative
coarse models do exist and may be accurate for all time.

\subsection{Compare with averaging}

\cite{Papavasiliou06},\S5, explored the multiscale, equation free,
modelling of the simple, two variable, one noise, stochastic system
\begin{eqnarray}
    dx&=&-(y+y^2)d\tau \,, \label{eq:pkx} \\
    dy&=&-\frac1\epsilon(y-x)d\tau +\frac1{\sqrt\epsilon}\,dW_\tau\,.
    \label{eq:pky}
\end{eqnarray}
This system has two time scales for small parameter~$\epsilon$: for
small~$\epsilon$ the fast variable~$y$ decays quickly to $y\approx x$
on a $\tau$~time scale~$\Ord{\epsilon}$; substituting this approximate
balance into~\eqref{eq:pkx} gives $dx\approx -(x+x^2)d\tau$ in the
absence of noise.
That is, over $\tau$~times longer than~$\Ord\epsilon$ the slow
variable~$x$ evolves.
We compare the information provided by averaging to that provided by
stochastic normal forms in multiscale modelling.

Many apply methods of singular perturbations to systems of the
form~\eqref{eq:pkx}--\eqref{eq:pky}. For example, \cite{Papavasiliou06}
use the method of averaging to deduce that
\begin{equation}
    x=\avx+\Ord{\sqrt\epsilon}
    \qtq{where}
    d\avx=-(\avx+\avx^2+\rat12)d\tau \,.
    \label{eq:pkavx}
\end{equation}
That is, solutions of~\eqref{eq:pkx}--\eqref{eq:pky} are modelled to an
error~$\Ord{\sqrt\epsilon}$ by the deterministic \ode~\eqref{eq:pkavx}
which applies over $\tau$~times longer than~$\Ord{\epsilon}$.
The noise in the fast variable~$y$ generates the extra
drift~$-\rat12\,d\tau$ in~\eqref{eq:pkavx} through the quadratic
nonlinearity in the slow equation~\eqref{eq:pkx}.
However, averaging gives no basis for improving the
$\Ord{\sqrt\epsilon}$~error: such errors are often large in
applications as the scale separation may only be an order of magnitude
or two; for example, \cite{Papavasiliou06} simulate
\sde{}s~\eqref{eq:pkx}--\eqref{eq:pky} with scale separation
$\epsilon=0.01$ implying errors are roughly~$\sqrt\epsilon=10\%$\,.
Nor does averaging recognise the stochastic fluctuations, seen in the Example \sde~\eqref{eq:sde0}, induced in the slow variable~$x$ through fluctuations in the fast variable~$y$.
Stochastic normal forms extract both effects, and more as well.

As in previous sections, computer algebra modified from that by
\cite{Roberts06j} readily derives a stochastic normal form for the
system~\eqref{eq:pkx}--\eqref{eq:pky}.  But first we avoid the
straightjacket of singular perturbations by simply rescaling time to
$t=\tau/\epsilon$\,: that is, we adopt a time scale~$t$ where the rapid
transients decay on a $t$~time of~$\Ord{1}$, and the slow variable~$x$
evolves on long times $\Delta t\sim 1/\epsilon$\,.  The example
system~\eqref{eq:pkx}--\eqref{eq:pky} is then identical to the
Stratonovich system
\begin{eqnarray}
    \dot x&=&-\epsilon(y+y^2) \,, \label{eq:pkxt} \\
    \dot y&=&-y+x +\sigma\phi(t)\,,
    \label{eq:pkyt}
\end{eqnarray}
when the new noise magnitude $\sigma=1$\,.  I introduce the noise
magnitude~$\sigma$ in the \sde\ system~\eqref{eq:pkxt}--\eqref{eq:pkyt} in
order to control truncation of noise effects.
Minor modifications of the previous computer algebra then
discovers that the stochastic coordinate transform
\begin{eqnarray}
    x&=& X
    +\epsilon(Y+\rat12Y^2+2XY)
    \nonumber\\&&{}
    +\epsilon\sigma \left[ (1+Y+2X)\Z-\phi +Y\Z+\phi\right]
    \nonumber\\&&{}
    +\rat12\epsilon\sigma^2(\Z-\phi)^2
    +\Ord{\epsilon^2}
    \,,\\
    y&=& Y+X +\sigma\Z-\phi
    -\rat12\epsilon Y^2
    \nonumber\\&&{}
    +\epsilon\sigma\left[ (1-Y+\Z-)\Z-\phi +Y\Z+\phi \right]
    \nonumber\\&&{}
    +\epsilon\sigma^2\Z-\left[\phi\Z-\phi+\rat12(\Z-\phi)^2 \right]
    +\Ord{\epsilon^2+X^2}
    \,,
\end{eqnarray}
maps the \sde\ system \eqref{eq:pkxt}--\eqref{eq:pkyt} into the following Stratonovich \sde\ system for the new variables $X$~and~$Y$:
\begin{eqnarray}
    \dot X&=& -\epsilon(X+\epsilon X+X^2)
    -\epsilon\sigma(1-2\epsilon+2X)\phi
    -\epsilon\sigma^2\phi\Z-\phi
    \nonumber\\&&{}
    +\Ord{\epsilon^4+X^4+\sigma^4} \,, \label{eq:pkxxd} \\
    \dot Y&=& Y\left[ (-1 +\epsilon +\epsilon^2 +2\epsilon^3)
    +(2\epsilon +4\epsilon^2)X
    +\sigma(2\epsilon+6\epsilon^2)\phi 
    \right]
    \nonumber\\&&{}
    +\Ord{\epsilon^4+X^4+\sigma^4} \,. \label{eq:pkyyd}
\end{eqnarray}
As before, the utility of this normal form transformation is that the
\sde~\eqref{eq:pkyyd} shows that the transformed fast variable~$Y\to0$
exponentially quickly from a wide range of initial conditions for small
scale separation~$\epsilon$.
Moreover, the methodology may refine the approximation, through further
iteration, to suit a wide range of specified finite scale
separation~$\epsilon$.

This normal form transformation also shows that the new slow
variable~$X$ evolves independently of the fast variable~$Y$,
see~\eqref{eq:pkxxd}, both throughout the initial transient as well as
thereafter: there are no initial transients in~$X$.
Furthermore, being just a transform form of the original
\sde~\eqref{eq:pkyt}, the \sde~\eqref{eq:pkxxd} for the slow
variable~$X(t)$ applies for all times, albeit to the truncation error;
in contrast, the averaged system generally only applies for a finite
time span.

Although not immediately apparent, the leading approximation of the slow~$X$ evolution~\eqref{eq:pkxxd}
is the averaged model~\eqref{eq:pkavx}.  The quadratic noise
term in~\eqref{eq:pkxxd} generates a mean drift and an effective new 
noise over long times: \cite{Roberts05c} and \cite{Chao95} argued that
over long times 
\begin{equation}
    \phi\Z-\phi \mapsto \rat12+\rat1{\sqrt2}\alt\phi \,,
    \label{eq:xlong}
\end{equation}
where $\alt\phi(t)$ is a new `white' noise process independent of the
original noise process~$\phi(t)$.  Thus the slow variable
\sde~\eqref{eq:pkxxd} is effectively the \sde
\begin{equation}
    \dot X= -\epsilon(\rat12\sigma^2+X+\epsilon X+X^2)
    -\epsilon\sigma(1-2\epsilon+2X)\phi
    -\epsilon\sigma^2\rat1{\sqrt2}\alt\phi \,.
    \label{eq:effxxd}
\end{equation}
Reverting to the original (slow) time~$\tau$, setting $\sigma=1$ to
match the original noise intensity, and re-expressing, the
\sde~\eqref{eq:effxxd} becomes
\begin{equation}
    dX=-(\rat12+X+X^2+\epsilon X)d\tau
    -\sqrt\epsilon(1-2\epsilon+2X)dW_\tau
    -\sqrt{\rat\epsilon2}\,d\alt W_\tau \,.
    \label{eq:effxd}
\end{equation}
The deterministic \ode\ found at leading order,
$dX=-(\rat12+X+X^2)d\tau$\,, is the averaged model~\eqref{eq:pkavx}.
However, the \sde~\eqref{eq:effxd} also makes explicit some of the
errors in averaging.
Firstly, the $\sqrt\epsilon$~error of the averaged
model~\eqref{eq:pkavx} comes from its neglect of the stochastic
fluctuations: to leading order we can combine noise processes
$W$~and~$\alt W$ to determine that the slow variables are better modelled
by the \sde\ $dX =-(\rat12+X+X^2)d\tau
+\sqrt{3\epsilon/2}\,d\aalt W_\tau$\,; although the two stochastic terms in the
\sde~\eqref{eq:effxd} are even better.
Secondly, the \sde~\eqref{eq:effxd} also discerns
$\Ord{\epsilon}$~contributions to the deterministic terms which may well have
significant effects at finite scale separation~$\epsilon$.
Simple averaging misses all of these effects.

\subsection{Equation free simulation}

\cite{Kevrekidis03b} promote a framework for computer aided, equation
free, multiscale analysis, which empowers systems specified at a
microscopic level of description to perform modeling tasks at a
macroscopic, systems level.  When the microscopic simulator is
stochastic, that is Monte--Carlo, or effectively stochastic, such as
molecular and discrete element simulators, then the issues addressed in
this article of the nature and extraction of long term dynamics from a
stochastic system are crucial to the equation free methodology.

Equation free modelling is designed to solve specific multiscale
systems with specific finite scale separations.  Thus a challenge for
future research is to maintain reasonable accuracy in estimating the
long term dynamics by extracting from numerical realisations the sort
of information extracted by these normal form coordinate transforms and
without knowing any algebraic representations of the systems of
interest.  The stochastic normal form transformation shows what might
be achieved in principle.  The challenge is to find out how to achieve it
from just a finite number of short bursts of realisations.

On the macroscale the stochastic effects may be relatively small.
However, a deterministic macroscale model is often structurally
unstable: one example is the structural instability of the averaged
model $d\bar x=0\,dt$ for the \sde~\eqref{eq:sde0}; instead we prefer the
stochastic model $dX=\epsilon\,dW$ of the \sde~\eqref{eq:Sde0}.
Moreover, even on the macroscale a deterministic model for some
averaged slow variable is almost inevitably different from the average
of the system with noise included.  This difference follows from the
same line of argument that establishes that the expectation of
realisations is generally different from the expected position of the
stochastic slow manifold, see~\eqref{eq:exy}.  Noise induced mean
drift must be recognised \cite[]{Srinamachchivaya90, Srinamachchivaya91}.

To do coarse step integration we need to estimate the macroscale drift
from the microscale simulations.  Noise hides this drift making
accurate estimation difficult.  We need either long bursts of
microsimulations, or many realisations, or variance reduction
techniques \cite[]{Papavasiliou06}, or a combination of all three.  

As well as the drift, the fluctuations in the macroscopic quantities
should be modelled.  Thus the macroscale integration should be that
of a system of \sde{}s.  Because the macroscale \sde{}s model
microscale processes, I conjecture that the macroscale \sde{}s must be
interpretted as Stratonovich \sde{}s.  The challenge is to develop
Stratonovich integrators that only use short bursts of realisations.

In equation free simulation one projects a macroscopic time step into
the future, then executes a burst of microscale simulation in order to
estimate the macroscopic rate of change \cite[]{Kevrekidis03b}.
Initial rapid transients must be ignored in each burst as the
microscopic system attains a quasi-equilibrium.  In a stochastic
system, the true \ssm\ can only be identified via integrals over fast
time scales, see Section~\ref{sec:aa}.  However, these are generally
integrals of both the past and the future.  Thus, to estimate
macroscopic rates of change in a stochastic system, we must not only
neglect initial transients, but also data from the end of a burst of
microscopic simulation in order to be able to account for the integrals
which anticipate the noise processes.

Lastly, the gap-tooth scheme empowers equation free modelling across
space scales as well as time scales \cite[e.g.]{Gear03}.  For
spatiotemporal stochastic systems we need theoretical support for the
notion that \spde{}s can be modelled by the gap-tooth scheme in the
same way as deterministic \pde{}s \cite[]{Roberts06d}.  Only then will
we be assured that we can cross space scales as well as time scales.


\section{Long time modelling of stochastic oscillations}
\label{sec:osc}

The previous sections focus on the separation of slow modes from
fast modes in stochastic systems.  Persistent oscillations are another
vitally important class of dynamics.  Hopf bifurcation is the
example considered in this section, but many other
cases occur including wave propagation.  The challenge addressed here
is how to consistently model the evolution of  oscillations
\emph{over long time scales} when the oscillations are fast and in the
presence of stochastic noise fluctuations over all time scales.  To model
over long time scales we  eliminate from the model \emph{all}
fast time dynamics.

As an example let us explore the stochastic Duffing--van der Pol
dynamics also analysed by \cite{Arnold93} and \cite{Arnold98}:
\begin{equation}
    \ddot x_1= (\alpha+\sigma\phi(t)) x_1+\beta \dot x_1
    -x_1^3-x_1^2\dot x_1\,,
    \label{eq:dvdp}
\end{equation}
where, as before, $\phi$~is some white noise process.  \cite*{Arnold95} describe
the importance of the stochastic system~\eqref{eq:dvdp} in applications.  
In the absence of noise, $\sigma=0$\,, this system exhibits
\begin{enumerate}
	\item a deterministic pitchfork bifurcation as the
	parameter~$\alpha$ crosses zero with $\beta$~fixed, say $\beta=-1$
	for definiteness; and

	\item a deterministic Hopf bifurcation as the parameter~$\beta$
	crosses zero with $\alpha$~fixed, say $\alpha=-1$ for definiteness.
\end{enumerate}
In the presence of noise, $\sigma>0$\,, computer algebra~\cite[]{Roberts06j} readily derives the stochastic normal form for the
Duffing--van der Pol equation~\eqref{eq:dvdp} near the stochastic pitchfork bifurcation when parameter~$\alpha$ crosses zero with fixed~$\beta$.
This section explores the issues arising when constructing a normal form for the Duffing--van der Pol equation~\eqref{eq:dvdp} near the stochastic Hopf bifurcation as the parameter~$\beta$ crosses zero with $\alpha=-1$\,.

\subsection{Approaches to stochastic Hopf bifurcation}
\label{sec:ashb}

\begin{figure}
    \centering
    \begin{tabular}{c@{}cc}
        & (a) $\beta=-0.1$ & (b) $\beta=+0.1$ \\
        \raisebox{20ex}{$\dot x_1$}&
        \includegraphics[width=0.45\linewidth]{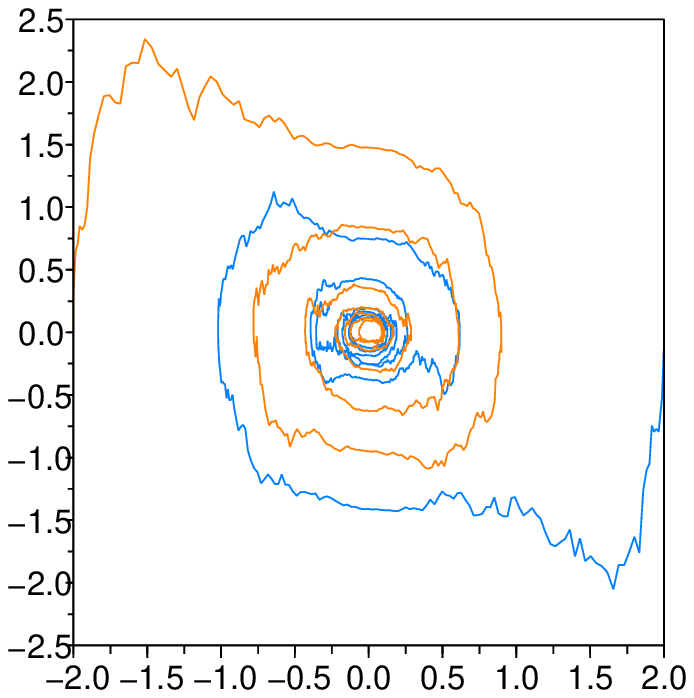} &
        \includegraphics[width=0.45\linewidth]{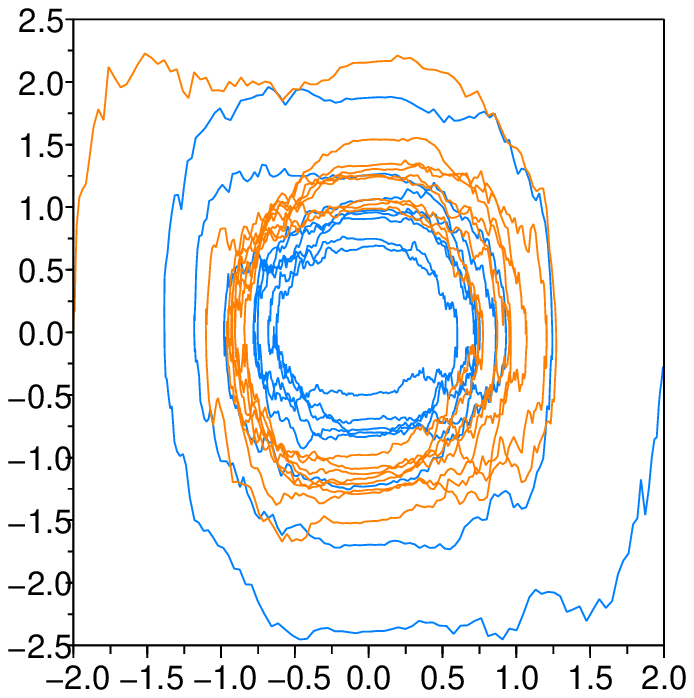} \\[-2ex]
        & $x_1$ & $x_1$
    \end{tabular}
	\caption{stochastic Hopf bifurcation in the Duffing--van der Pol
	oscillator~\eqref{eq:dvdp} as parameter~$\beta$ crosses zero with
	$\alpha=-1$ and noise amplitude $\sigma=0.5$\,.   Two
    realisations are plotted in each case.}
    \label{fig:shopf}
\end{figure}

Figure~\ref{fig:shopf} shows the noisy Duffing--van der Pol
oscillator~\eqref{eq:dvdp} of \cite{Arnold98} for parameter
$\beta=\pm0.1$\,.
Deterministically ($\sigma=0$), as parameter~$\beta$ crosses zero, a
Hopf bifurcation occurs from the stable equilibrium at the origin to a
stable limit cycle with frequency~$1$.
What happens in the presence of parametric stochastic forcing when
$\sigma\neq0$\,? Figure~\ref{fig:shopf} reaffirms that a noisy version
of the stochastic bifurcation takes place.
\cite{Coullet83} first explored a normal form of Hopf bifurcations with
noise.
Further research on such stochastic bifurcations elucidated some
fascinating fine structure.
For example, \cite{Keller99} explored the structure of the random
`limit cycle' attractor using a stochastic version of the subdivision
algorithm of \cite{Dellnitz97}; whereas \cite{Arnold98} explored the
structures using a normal form approach very close to that used here.
However, the emphasis here is not on the stochastic Hopf bifurcation as
such, but instead using it as the simplest prototype system with
stochastic oscillatory dynamics.
We look at the issues afresh to \emph{explore the characteristics of a
long term stochastic model} of such stochastic oscillatory dynamics.
In the future, these considerations will underpin the multiscale
modelling of stochastic oscillations and waves.

Section~\ref{sec:csnf} constructs a stochastic coordinate transform
from which we may extract significant properties of a stochastic
Hopf bifurcation.  Solutions of the Duffing--van der Pol
oscillator~\eqref{eq:dvdp} are most conveniently represented in complex exponentials as
\begin{equation}
    x_1\approx a(t)e^{it}+b(t)e^{-it}\,,
    \label{eq:oscab}
\end{equation}
where for real solutions~$x_1$, the amplitudes $a$~and~$b$ are complex
conjugates.  Then Section~\ref{sec:csnf} finds a stochastically
forced Landau model governing the evolution of the complex amplitudes
$a$~and~$b$:
\begin{eqnarray}
    \dot a&\approx&\rat12\beta a -(\rat12-\rat32i)a^2b
    +\sigma\sqrt{\delta/2}(a\phi_0-b\phi_{+2}) \,,
    \label{eq:ha}
    \\
    \dot b&\approx&\rat12\beta b -(\rat12+\rat32i)ab^2
    +\sigma\sqrt{\delta/2}(b\phi_0-a\phi_{-2}) \,,
    \label{eq:hb}
\end{eqnarray}
to errors $\Ord{\beta^2+\sigma^2+\epsilon^4}$
where $\epsilon=|(a,b)|$ measures the size of the oscillations, and
where $\phi_m(t)$ are independent `white' noises arising from the
forcing components near frequencies $0$~and~$\pm2$ in the applied noise
process~$\phi(t)$; `near' means within~$\pm\delta$ of the specified
frequency.  This model resolves the slow evolution of the complex
amplitudes near the Hopf bifurcation, small~$\beta$, under the
influence of the nonlinearity and a weak stochastic forcing,
small~$\sigma$.  This model empowers long term simulations with
efficient large time steps as the complex amplitudes are
slowly-varying.

Note: the analysis also applies in the case of the deterministic
forcing~$\phi=\cos2t$\,, for which $\phi_0=0$ and $\sqrt{2\delta}
\phi_{\pm2}=\rat12$\,.  Then the above model, $\dot a \approx \rat12
\beta a +\rat14 \sigma b$ and its conjugate, successfully predicts
the Mathieu-like instability with eigenvalues $\lambda =\rat12\beta
\pm\rat14\sigma$\,.

View~\eqref{eq:oscab} as a time dependent coordinate transform of the
$(x_1,\dot x_1)$~phase plane.  In principle, any dynamics in the phase
plane may be described by the evolution of the complex amplitudes
$a$~and~$b$.  The utility of the coordinate transform~\eqref{eq:oscab}
is that it empowers a simple description of oscillations with frequency
near~$1$: namely \eqref{eq:ha}--\eqref{eq:hb} for the Duffing--van der
Pol oscillator~\eqref{eq:dvdp}.  However, to simply describe such
nonlinear stochastic oscillations Section~\ref{sec:csnf} modifies the
coordinate transform~\eqref{eq:oscab} through nonlinear and stochastic
terms.  That is, there is a time dependent, coordinate transform of the
phase plane that leads to the normal form \eqref{eq:ha}--\eqref{eq:hb}.

I emphasise this different view of~\eqref{eq:oscab}.  Many would view~\eqref{eq:oscab} as an approximation to~$x(t)$ that can only resolve slowly varying oscillations.  In contrast, I present~\eqref{eq:oscab} as the leading term in a coordinate transform, a reparametrisation, of the entire phase $(x_1,\dot x_1)$~plane that in principle encompasses \emph{all} dynamics in the phase plane.  The approximate model then arises by finding parameter regimes, in this new coordinate system, where the evolution of `coordinates' $a$~and~$b$ is usefully slow.

\paragraph{Amplitude/phase models do not decouple}
\cite{Arnold98} analysed a Hopf bifurcation by transforming to real
amplitude~$r$ and phase angle~$\varphi$ coordinates and deducing a
model $\dot r=\cdots$ and $\dot\varphi=1+\cdots$\,.
This approach is certainly effective for unforced deterministic
problems~\cite[e.g.]{Roberts06a}.
However, the presence of time dependent forcing, whether stochastic or
deterministic, breaks time translation symmetry.
Consequently, \cite{Arnold98} must couple the phase~$\varphi$ back into
the amplitude~$r$ evolution, as also seen in the normal form~(39) of
\cite{Arnold95}.
Such coupling of the fast phase into the notionally slow amplitude
confounds our aim to use the normal form for long time modelling.

Because of their different aim, \cite{Arnold95} convert back to a pair
of fast Cartesian variables to obtain a canonical system that is generic for
the class of stochastic Hopf bifurcations; thus they establish that the
pattern of behaviour they explore is generic for Hopf bifurcations.
But our aim here is different: we aim to construct models suitable for
exploring long time evolution; our normal form is consequently
different.
We use complex amplitude coordinates, the $a$~and~$b$ seen in
\eqref{eq:ha}~and~\eqref{eq:hb}, as originally proposed by
\cite{Coullet83}.


Stochastic averaging seems to suffer the same defect of not recognising
the broken time symmetry \cite[equations~(16--20)]{Arnold95}.
Stochastic averaging also does not appear to detect the split in Lyapunov
exponents present in stochastic Hopf bifurcations.

\paragraph{Prefer a strong model}
\cite{Olarrea96} comment that ``When the reduction to the normal form
is done \ldots\ only the deterministic part of the equations retain the
characteristic radial symmetry.'' and then assert ``This makes it
necessary to work with the two-dimensional probability distribution.''
Thus they introduce early in their analysis some probability
distributions governed by Fokker--Planck equations and hence derive only weak models.
In contrast, here we maintain strong modelling of each realisation of the
noise.
We avoid weak models.

\subsection{Construct a stochastic normal form}
\label{sec:csnf}

To construct the stochastic normal form for the Duffing--van der Pol
oscillator~\eqref{eq:dvdp}, with parameter $\alpha=-1$\,, I use an
iterative scheme to construct a useful nonlinear coordinate transform.
The coordinate transform must be time dependent to adapt to both the
oscillations and to the stochastic effects.
The starting approximation to the linear time dependent coordinate
transform is~\eqref{eq:oscab}.
Iterative modifications to~\eqref{eq:oscab} result in a description of
the Duffing--van der Pol oscillator~\eqref{eq:dvdp} which only has slow
processes suitable for long time simulation.

\paragraph{The homological equation}
Each step in the iteration improves the normal form description of
the dynamics.  Suppose that at some step in the iteration, the coordinate 
transform and consequent evolution is
\begin{displaymath}
    x_1=\xi(a,b,t)
    \qtq{where}
    \dot a=g(a,b,t)
    \qtq{and}
    \dot b=h(a,b,t)
\end{displaymath}
for some known functions~$\xi$, $g$~and~$h$.  Seek small corrections,
denoted by dashes, to~$\xi$, $g$~and~$h$ so that
\begin{equation}
    x_1=\xi+\xi'(a,b,t)
    \qtq{where}
    \dot a=g+g'(a,b,t)
    \qtq{and}
    \dot b=h+h'(a,b,t)
    \label{eq:oscup}
\end{equation}
better satisfies the Duffing--van der Pol oscillator~\eqref{eq:dvdp}.
We measure how well the Duffing--van der Pol oscillator~\eqref{eq:dvdp}
is satisfied by its residual,~$\res_{\protect\ref{eq:dvdp}}$.
Substitute~\eqref{eq:oscup} into the Duffing--van der Pol
oscillator~\eqref{eq:dvdp}, omit products of small corrections,
approximate $\xi \approx ae^{it}+be^{-it}$ and $g \approx h \approx
\beta \approx \sigma \approx 0$ whenever multiplied by a correction,
and deduce that in the complex amplitude coordinates, the homological
equation is
\begin{displaymath}
    \xi'_{tt}+\xi' +(i2g'+g'_t)e^{it} +(-i2h'+h'_t)e^{-it}
    +\res_{\protect\ref{eq:dvdp}}=0\,.
\end{displaymath}
But there is one further refinement: we aim for $\dot a=g$ and $\dot
b=h$ to only possess \emph{slow} dynamics; thus, presuming this aim is
possible, also omit the time derivatives $g'_t$~and~$h'_t$ to give the
homological equation
\begin{equation}
    \xi'_{tt}+\xi' +i2g'e^{it} -i2h'e^{-it}
    +\res_{\protect\ref{eq:dvdp}}=0\,.
    \label{eq:oschom}
\end{equation}
This approach avoids difficulties that appear in the homological
equation for amplitude-phase models.  The homological
equation~\eqref{eq:oschom} governs corrections to the complex
coordinate transform.

\subsection{Linear noise effects}

An iterative scheme to find a stochastic coordinate transform and
corresponding evolution was coded into computer algebra
\cite[]{Roberts06j}.  Iterative improvements to the coordinate
transform and the model continue until the residual of the Duffing--van
der Pol oscillator~\eqref{eq:dvdp} reaches a specified order
of error.  To effects linear in the noise magnitude~$\sigma$ the
iteration finds the stochastic model
\eqref{eq:ha}~and~\eqref{eq:hb} to the specified errors.  In terms of
the Fourier transform~$\tilde\phi(\freq)$ of the noise,
$\phi(t)=\int_{-\infty}^\infty e^{i\freq t} \tilde\phi(\freq)
\,d\freq$\,, the
corresponding stochastic complex coordinate transform is
\begin{eqnarray}
    x_1&=&
    ae^{it}+be^{-it}
    +\rat18\big[(1+i)a^3e^{i3t} +(1-i)b^3e^{-i3t} \big]
    \nonumber\\&&{}
    -\sigma ia\xint \frac1{\freq(\freq+2)} e^{i(\freq+1) t}
    \tilde\phi(\freq) \,d\freq
    \nonumber\\&&{}
    +\sigma ib\xint \frac1{\freq(\freq-2)} e^{i(\freq-1) t}
    \tilde\phi(\freq) \,d\freq
    \nonumber\\&&{}
    +\sqrt{2\delta} \sigma\Big[ \rat i4(a\phi_0-b\phi_2)e^{it}
    -\rat i4(b\phi_0-a\phi_{-2})e^{-it}
    \nonumber\\&&\quad{}
    -\rat i8a\phi_2e^{i3t}
    +\rat i8b\phi_{-2}e^{-i3t}
    \Big]
    +\Ord{\beta^2+\sigma^2+\epsilon^4,\delta^{3/2}}\,,
    \label{eq:x}
\end{eqnarray}
where the integration domain~$D$ avoids singularities in the integrand as explained in Section~\ref{sec:nrf}.

\subsubsection{Deterministic effects}

The first line of~\eqref{eq:x} describes the well established
deterministic shape of the limit cycle in the deterministic Hopf
bifurcation.
When the residual~$\res_{\protect\ref{eq:dvdp}}$ has terms with factors
$e^{imt}$ for some integer~$m$, $|m|\neq1$\,, and no other explicit
time dependence, then as usual we update the complex coordinate
transform by a correction~$\xi'$ proportional to $e^{imt}/(m^2-1)$, and
do not change the evolution, $g'=h'=0$\,.

Deterministic terms in the residual with factors~$e^{\pm it}$, and no
other explicit time dependence, such as the term $(i\beta
a+3a^2b)e^{it}$, are resonant and as usual must be assigned to correct
the evolution, upon dividing by the $\pm2i$~factor of the homological
equation~\eqref{eq:oschom}; see the deterministic nonlinear and
$\beta$~terms in the model~\eqref{eq:ha}--\eqref{eq:hb}.

\subsubsection{Non-resonant fluctuations}
\label{sec:nrf}

The second two lines of the transform~\eqref{eq:x} describe how
stochastic fluctuations non-resonantly perturb the oscillating
dynamics.  These arise from terms in the
residual~$\res_{\protect\ref{eq:dvdp}}$ of the form
\begin{displaymath}
    \phi(t)e^{\pm it}=\int_{-\infty}^\infty e^{i(\freq\pm1)t}
    \tilde\phi(\freq) \,d\freq\,.
\end{displaymath}
Away from resonance, namely in the domain $D=\mathbb R\backslash
\cup_{m\in\{-2,0,2\}}[m-\delta,m+\delta]$\,, these terms in the
residual generate the desingularised integrals
in~\eqref{eq:x}.\footnote{In analyses to higher order in the
oscillation amplitude more resonant frequencies occur; for example,
integrals arise with singularities at frequencies $\freq=\pm4$ in some
terms of~$\Ord{\epsilon^2\sigma}$.
In such higher order analyses the domain of integration~$D$ will have
further intervals excised to avoid resonances.} Rewriting these
integrals as a convolution~$f(t)\star\phi(t)$ recognise that formally
$f=e^{\pm it}\xint \frac1{\freq(\freq\pm 2)} e^{i\freq t} \,d\freq$\,.
This integral for the convolution kernel~$f$ may be written in terms of
the Sine integral \cite[\S5.2]{Abramowitz64} from which we deduce that
the convolution kernel~$f(t)$ decays like~$1/(\delta|t|)$ for
large~$|t|$.
Assuming that convolutions of~$f(t)$ with stochastic white noise do
converge in some sense, the complex transform appears to necessarily involves the entire past and future of the noise.
In contrast to the pitchfork bifurcation, which only needs to look a
little way into the future and the past, in the Hopf bifurcation we
look far into the future and the past in order to construct the
stochastic coordinate transform.
    
In contrast, \cite{Coullet83}, in their equations (18)~and~(19), assign
the entire integral to the evolution~\eqref{eq:ha}--\eqref{eq:hb},
rather than to the transformation, just because one frequency is
resonant.
This approach seems inconsistent in the neglect of the time derivatives
$g'_t$~and~$h'_t$ in the homological equation as such derivatives are
large for `white' noise.
Their assignment to the evolution is consistent when the noise
$\phi(t)$ has a narrow band spectrum around the resonant frequencies.

\begin{figure}
    \centering
    \begin{tabular}{c@{}c}
    \raisebox{20ex}{$\phi_0(t)$} &
    \includegraphics{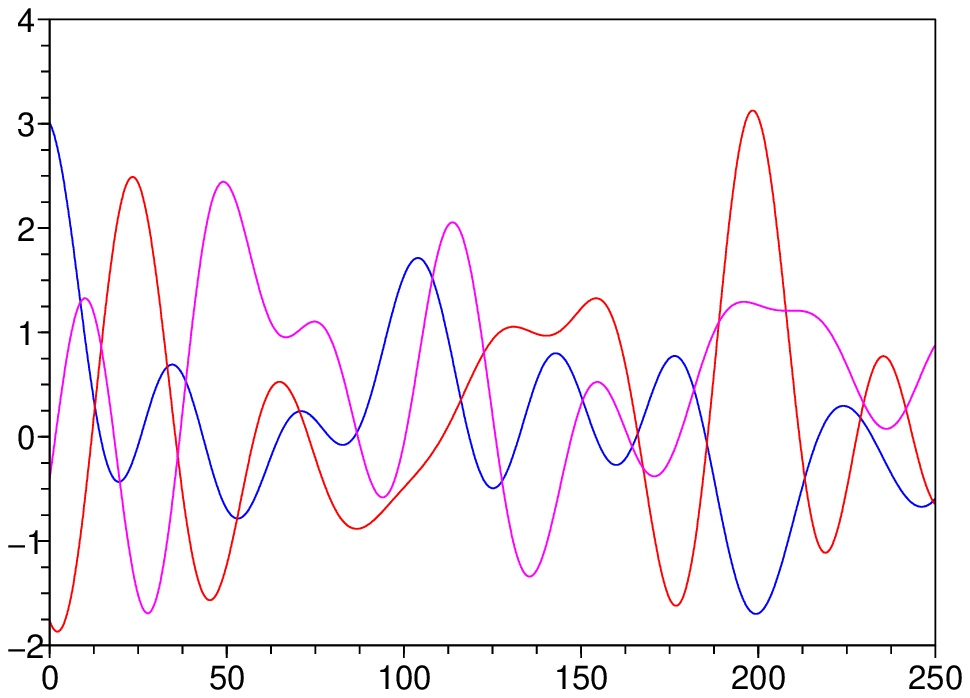} \\[-2ex]
    & $t$
    \end{tabular}
	\caption{schematic plot of three realisations of the amplitude of
	resonant noise~$\phi_0(t)$ for some mesoscale cutoff~$\delta$.}
    \label{fig:meso}
\end{figure}

\subsubsection{Resonant fluctuations}
\label{sec:rf}

The excised parts of the integrals in the transform~\eqref{eq:x}
correspond to resonances.  These resonances generate terms in the
model~\eqref{eq:ha}--\eqref{eq:hb} involving components of the
(complex) noise process
\begin{equation}
	 \phi_m(t)=\frac1{\sqrt{2\delta}} \int_{m-\delta}^{m+\delta}
	 e^{i(\freq-m)t} \tilde\phi(\freq) \,d\freq\,,
\end{equation}
normalised so that $\E{|\phi_m|^2}=1$ under the original white noise
assumption that $\E{ \cc{\tilde\phi(\freq)} \tilde\phi(\alt\freq)}
=\delta(\freq-\alt\freq)$ (here $\delta()$~denotes the Dirac delta function and $*$~the complex conjugate);
Figure~\ref{fig:meso} plots three realisations.
Being a narrow band integral (with the dominant frequency accounted for
by the $e^{-imt}$~factor) the $\phi_m(t)$ are slowly varying noise
processes: Figure~\ref{fig:meso} shows $\phi_0(t)$ has slow variations on
the fast times scale $\Delta t=2\pi$ of the oscillations.
They are independent of each other as the domains of integration do not
overlap (for small cutoff~$\delta$).
Each $\phi_m(t)$ has autocorrelations which decay on a time scale of
order~$1/\delta$, roughly the width of the window in
Figure~\ref{fig:meso}, but for time scales${}\gg1/\delta$ the
autocorrelation is zero and the $\phi_m$ look like white noise
processes.
Thus choose the `cutoff'~$1/\delta$ to be a mesoscopic time scale: one
longer than the period of the limit cycle; but much shorter than the
long macroscopic time scale on which the
model~\eqref{eq:ha}--\eqref{eq:hb} is to be used.
Then $\phi_m(t)$ are effectively independent white noise processes in
the long term model.
 
Encouragingly, although the Fourier transform~$\tilde\phi(\freq)$
requires the entire history of the noise, the parts of it that appear
in the model~\eqref{eq:ha}--\eqref{eq:hb} are essentially local in time.
That is, as for non oscillatory dynamics, \emph{the long term model
itself does not require anticipation of the noise.}

The fourth and fifth lines in the transform~\eqref{eq:x} arise through
the excision of the resonant parts of the frequency domain from the
integrals in the coordinate transform~\eqref{eq:x}.
These resonant frequencies not only affect the evolution but also the
coordinate transform as seen in these two lines of~\eqref{eq:x}.

These resonant fluctuations also force the complex amplitudes
$a$~and~$b$ to change their meaning in the presence of noise.
I do not precisely and explicitly define the complex amplitudes
$a$~and~$b$; \emph{implicitly} they are the component in $e^{\pm it}$
in the oscillations.
However, whatever definition one may try to adopt, implicitly or
explicitly, the noise changes the definition through the terms
appearing on the fourth and fifth lines in the transform~\eqref{eq:x}.
Recall that in non-oscillatory systems noise also changes the presumed
definition of slow variables: for two examples, the \ssm{}s
\eqref{eq:sde0nf} and~\eqref{eq:toyssm}
show that we cannot parametrise a \ssm\ in terms of the original
slow variable~$x$, but a new variable~$X$ which is necessarily
different in the presence of noise.
Similarly here: in the presence of noise, the coefficient of~$e^{it}$
in the stochastic coordinate transform is not just the complex
amplitude~$a$ but instead is approximately~$a+i\sigma \sqrt{2\delta}
\rat14 (a\phi_0 -b\phi_2)$\,, and analogously for the coefficient
of~$e^{-it}$.
Noise affects the meaning of the complex amplitudes.

These terms of the fourth and fifth lines in the
transform~\eqref{eq:x}, and the corresponding terms in the
model~\eqref{eq:ha}--\eqref{eq:hb}, are proportional to~$\sqrt{\delta}$
where $\delta$~is the small width of the domain excised from frequency
space about the resonant terms.
Can these terms be ignored as small? I contend it depends upon the use
of the slow model~\eqref{eq:ha}--\eqref{eq:hb}.
In a long term simulation we may use macroscopic time steps of
size~$\Delta t$, say, in numerically
integrating~\eqref{eq:ha}--\eqref{eq:hb}.
In this numerical integration we would treat the $\phi_m(t)$~noises as
white; thus their decorrelation time~$1/\delta$ must be less than the
numerical time step~$\Delta t$.
That is, a lower bound for the excised mesoscale cutoff is
$\delta>1/\Delta t$\,.
Thus, a stochastic time integrator could treat these terms as
of~$\Ord{1/\sqrt{\Delta t}}$ but no smaller.

%

\subsection{Quadratic noise effects}

In many applications, quadratic noise effects generate important mean
deterministic drifts \cite[e.g.]{Srinamachchivaya90,
Srinamachchivaya91}.
This is easily seen in some examples, even using the method of
averaging \cite[\S5, e.g.]{Papavasiliou06}.
Such mean drifts are often important.
Thus, generically we must also explore how to analyse quadratic noise
effects.

\subsubsection{Double integrals of noise complicate}

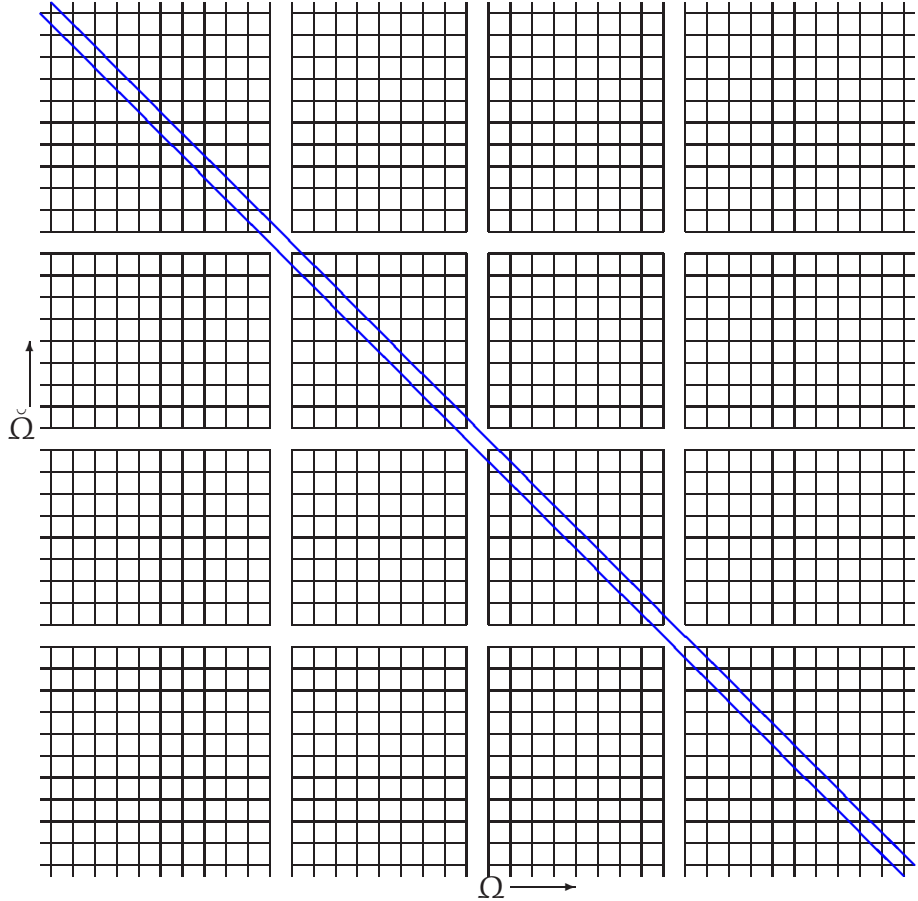
\begin{figure}
\centering
    \setlength{\unitlength}{0.8ex}
\begin{picture}(80,80)(-40,-40)
    \put(0,-42){$\freq$}
    \put(3,-41){\vector(1,0){6}}
    \put(-43,0){$\alt\freq$}
    \put(-41,3){\vector(0,1){6}}
    \multiput(-39,-40)(2,0){40}{
    \put(0,0){\line(0,1){21}}
    \put(0,23){\line(0,1){16}}
    \put(0,41){\line(0,1){16}}
    \put(0,59){\line(0,1){21}}
    }
    \multiput(-40,-39)(0,2){40}{
    \put(0,0){\line(1,0){21}}
    \put(23,0){\line(1,0){16}}
    \put(41,0){\line(1,0){16}}
    \put(59,0){\line(1,0){21}}
    }
    \color{blue}\thicklines
    \put(-40,39){\line(1,-1){79}}
    \put(-39,40){\line(1,-1){79}}
\end{picture}
    \caption{the integration domain $D\times D$, hatched,  also
    has a further resonant region, the diagonal blue strip, excised to
    give the integration domain~$\alt D$ for double integrals over the
    noise frequency.}
    \label{fig:domdom}
\end{figure}
For oscillatory dynamics, as in the Hopf bifurcation of the
Duffing--van der Pol oscillator~\eqref{eq:dvdp}, the outstanding
complication is the appearance of double integrals across all
frequencies in the stochastic fluctuations.
Quadratic noise effects \emph{not} involving such double integrals are
straightforwardly handled as before.
Terms of~$\Ord{\sigma^2}$ will contain double integrals of the form
$\int_D \int_D \cdot d\freq \,d\alt\freq$ where both
$\freq$~and~$\alt\freq$ represent noise frequencies.
The (black) hatched region in Figure~\ref{fig:domdom} shows this domain
of integration.
However, in the Hopf bifurcation of the Duffing--van der Pol
oscillator~\eqref{eq:dvdp}, the kernel of such double integrals also
has a singularity along the line $\freq+\alt\freq=0$\,.
Thus excise the (blue) diagonal strip shown in Figure~\ref{fig:domdom}
to remove the singularity to leave an integral over non-resonant
effects in the domain~$\alt D$.
Then additionally analyse the excised strip as a resonant effect that
directly influence the evolution of complex amplitudes $a$~and~$b$.

Recall we use the residual of an \sde\ system to drive corrections to the normal form stochastic coordinate transform.
In the residual of the Duffing--van der Pol oscillator~\eqref{eq:dvdp} quadratic noise terms arise of the form 
\begin{displaymath}
     \int_D\int_D e^{i(\freq+\alt\freq\pm1)t}
    K_\pm(\freq,\alt\freq) \tilde\phi(\freq)
    \tilde\phi(\alt\freq)\,d\freq\,d\alt\freq\,,
\end{displaymath}
where the integrand kernels are
\begin{equation}
    K_\pm = -\frac{(\freq+\alt\freq \pm \freq\alt\freq)(\freq+\alt\freq\pm2)}
    {2(\freq\pm2)(\alt\freq\pm2)\freq\alt\freq} \,.
\end{equation}
Before excising the blue strip in Figure~\ref{fig:domdom} to avoid the division by zero near $\freq+\alt\freq=0$\,, change the parametrisation of the integration domain to $\ff
=\rat12(\freq-\alt\freq)$ and $\alt\ff =\freq+\alt\freq$ so that $\freq =\ff
+\rat12\alt\ff$\,, $\alt\freq =-\ff +\rat12\alt\ff$\,, and the Jacobian of the
transform is one: parameter~$\alt\ff$ measures the distance from resonance.   In this new parametrisation, the integration kernels
\begin{eqnarray}
    K_\pm&=&
    \frac{2(4\ff^2\mp4\alt\ff-{\alt\ff}^2)(2\pm\alt\ff)}
    {(2\ff\pm4+\alt\ff)(2\ff\mp4-\alt\ff)(2\ff+\alt\ff)(2\ff-\alt\ff)}
    \\&\to& \frac{1}{(\ff+2)(\ff-2)}
    \quad\text{as }\alt\ff\to0\,.
    \nonumber
\end{eqnarray}
Then the double integrals in the residual are split into non-resonant and resonant parts:
\begin{eqnarray}
    I_\pm&=&\iint_{\alt D} e^{i(\freq+\alt\freq\pm1)t}
    K_\pm(\freq,\alt\freq) \tilde\phi(\freq)
    \tilde\phi(\alt\freq) \,d\freq\,d\alt\freq
    \nonumber\\&&{}
    +e^{\pm it} \int_{-\delta}^\delta e^{i\alt\ff t}\tilde\psi_\pm(\alt\ff) \,d\alt\ff \,,
    \label{eq:iik}
\end{eqnarray}
where
\begin{equation}
    \tilde\psi_\pm(\alt\ff)=\int_D 
    K_\pm\,
    \tilde\phi(\ff+\rat{\alt\ff}2) 
    \tilde\phi(-\ff+\rat{\alt\ff}2) \,d\ff\,,
    \label{eq:ipsi}
\end{equation}
and where domain~$\alt D=D\times D$ without the resonant strip as excised 
in Figure~\ref{fig:domdom}.  
The non-resonant double integral in~\eqref{eq:iik} contributes components to the stochastic coordinate transform.  
The resonant integral on the second line of~\eqref{eq:iik} contributes a component to the evolution in the new coordinates. 
Although the details will differ, the above integrals will appear in the analysis of general stochastic Hopf bifurcations.

\begin{figure}
    \centering
    \begin{tabular}{c@{}c}
       \raisebox{20ex}{$\psi_{r,i}$} &
       \includegraphics{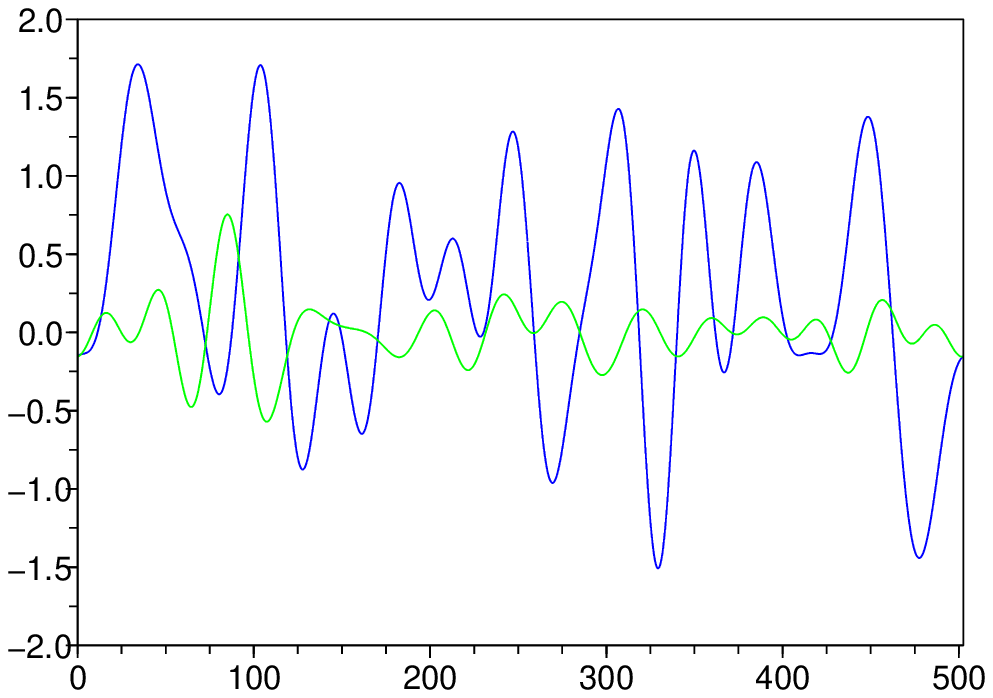}  \\[-2ex]
       & $t$
    \end{tabular}
	\caption{one realisation of the complex quadratically generated
	`noise' $\psi_\pm(t)\approx 0.87\psi_r(t)\pm i0.20\psi_i(t)$ where
	the real part is the larger blue curve and the imaginary part is
	the smaller green curve.  The resonant window size $\delta=0.2$\,.}
    \label{fig:cpvddt}
\end{figure}

The stochastic dynamics in the normal form coordinates will involve the integral~\eqref{eq:ipsi}.
The integral~\eqref{eq:ipsi} specifies the Fourier transforms of two complex conjugate components~$\psi_\pm(t)$ that express a nonlinear combination of the original noise process~$\phi(t)$.  
Here write these in terms of the real and imaginary parts 
\begin{equation}
    \psi_\pm(t)=c_r\psi_r(t)\pm i c_i\psi_i(t)\,,
\end{equation}
where the constants $c_r$~and~$c_i$ are chosen so the variances $\E{\psi_r^2}=\E{\psi_i^2}=1$\,; these constants do not seem to vary significantly with mesoscale cutoff~$\delta$.
Figure~\ref{fig:cpvddt} shows one realisation of~$\psi_\pm(t)$ illustrating that~$\psi_\pm(t)$ vary slowly over one period of the microscale limit cycle, and that they look like white noise processes over the long time scales resolved by the complex amplitudes $a$~and~$b$.  
In the Hopf bifurcation of the Duffing--van der Pol oscillator~\eqref{eq:dvdp} the processes $\psi_\pm(t)$ appear to have zero mean; this may not hold for other stochastic Hopf bifurcations.

\subsubsection{Refine the normal form transformation}

Separating the double integrals as described, computer algebra
\cite[]{Roberts06j} iteratively refines the stochastic coordinate
transform~\eqref{eq:x} and simultaneously derives the following \sde{}s
for the evolution of the complex amplitudes of the Duffing--van der Pol
oscillator~\eqref{eq:dvdp}:
\begin{eqnarray}
    \dot a&\approx&\rat12\beta a -(\rat12-\rat32i)a^2b
    +\sigma\sqrt{\delta/2}(a\phi_0-b\phi_{+2}) 
    \nonumber \\&&{}
    +i\rat12\sigma^2(c_r\psi_r+ic_i\psi_i)a
    -i\delta\sigma^2(\rat14\phi_0^2+\rat18\phi_2\phi_{-2})a
    \,,  \label{eq:ha2}
    \\
    \dot b&\approx&\rat12\beta b -(\rat12+\rat32i)ab^2
    +\sigma\sqrt{\delta/2}(b\phi_0-a\phi_{-2}) 
    \nonumber \\&&{}
    -i\rat12\sigma^2(c_r\psi_r-ic_i\psi_i)b
    +i\delta\sigma^2(\rat14\phi_0^2+\rat18\phi_2\phi_{-2})b
    \,. \label{eq:hb2}
\end{eqnarray}
The order of error in these \sde{}s
is~$\Ord{\epsilon^4+\sigma^3+\beta^2,\delta^{3/2}}$.
These \sde{}s account for more noise interactions than the lower order
model \eqref{eq:ha}--\eqref{eq:hb} and thus are more accurate.

For very small mesoscale cutoff~$\delta$, that is for simulations on
very long time scales, the quadratic noise effects involving
$\psi_r$~and~$\psi_i$ are the dominant influences on the complex
amplitudes $a$~and~$b$ of the oscillations of the Duffing--van der Pol
oscillator~\eqref{eq:dvdp}.
These two noise processes, see the integral~\eqref{eq:ipsi}, arise as
integrals of quadratic terms in the original noise process~$\phi$.
Analogously to the quadratic noise processes analysed on stochastic
slow manifolds \cite[\S5]{Roberts05c}, as used in~\eqref{eq:xlong}, I
conjecture that $\psi_r$, $\psi_i$ and~$\phi$ are effectively
independent when sampled over long time scales.
Consequently, over very long time scales, one would model real dynamics
of the stochastic Duffing--van der Pol oscillator~\eqref{eq:dvdp} by
the Stratonovich \sde
\begin{equation}
    da\approx \big[ \rat12\beta a -(\rat12-\rat32i)|a|^2a \big] dt
    +i\rat12c_r\sigma^2a\,dW_r -\rat12c_i\sigma^2a\,dW_i\,,
\end{equation}
where complex~$a$ measures the amplitude and phase of the oscillations,
$W_r$~and~$W_i$ denote independent Wiener processes, and $c_r\approx
.87$ and $c_i\approx .20$ (Figure~\ref{fig:cpvddt}).

For medium mesoscale cutoff~$\delta$ use the more complete \sde\
model~\eqref{eq:ha2}.
This model, with its effects in $\sqrt\delta$~and~$\delta$, will be
needed when the desired time resolution of a numerical simulator,
essentially the integrator's time step~$\Delta t$, is within a few
orders of magnitude of the natural period of oscillations, here the
period is about~$2\pi$.
A challenge for future research is to construct special \sde\ numerical
iteration schemes when, as here, \emph{the \sde\ itself depends upon
the chosen time step~$\Delta t$};
I am only aware of \sde\ schemes which assume the \sde\ is independent
of the time step \cite[e.g.]{Higham05, Kloeden92}.
Physically, the dependence upon the macroscopic time step is due to
the difficulty in discerning what is and is not a resonant forcing of
the oscillations, see Sections~\ref{sec:nrf}--\ref{sec:rf}.
In multiscale modelling, as shown here, the macroscopic system, whether
expressed as algebraic equations or solved using equation free methods
\cite[]{Kevrekidis03b}, may depend upon the the length or time scale
chosen for simulation.

The specific equations and formulae in the section are specific to the
Duffing--van der Pol equation~\eqref{eq:dvdp}.  Nonetheless, I contend
that the nonlinear and stochastic nature of these Duffing--van der Pol
oscillations are generic for most of the interesting issues discussed
in this section.  Consequently, I conjecture that almost all long time
scale modelling of stochastic oscillations has to address and resolve the
issues discussed in this section.

\section{Conclusion}

Stochastic coordinate transforms illuminate the modelling of multiscale stochastic systems.  
Being a coordinate transform, a resultant `stochastic normal form' describes the complete dynamics of the original system, Proposition~\ref{thm:nf}.  
From the normal form we easily extract the stochastic slow dynamics that are of interest over macroscopic times, from the uninteresting fast dynamics \cite[\S8.4, e.g.]{Arnold03}.  This approach is more powerful than averaging as the coordinate transform may be systematically refined, especially with the aid of computer algebra \cite[]{Roberts06j}, and so errors are more controlled.

In contrast to earlier work, this article argues that two modelling simplifications may always be achieved without sacrificing fidelity with the original stochastic system.  
Firstly, the stochastic slow manifold and the evolution thereon need not have any terms anticipating the original noise processes, Proposition~\ref{thm:mem}.  
Secondly, effects linear in the noise processes in the evolution on the stochastic slow manifold need not involve any memory integrals either,  Proposition~\ref{thm:nf}. 
Section~\ref{sec:eg} illustrates these principles for the example \sde\ system \eqref{eq:toy0}.

A challenge for future research is to let the algebraic techniques used herein inspire development of numerical techniques useful for multiscale computations.  
From a finite number of bursts of stochastic realisations we need to determine information to empower making macroscale time steps while remaining faithful to the underlying stochastic dynamics. 

Section~\ref{sec:osc} explored oscillatory dynamics in the stochastic Duffing--van der Pol equation~\eqref{eq:dvdp}.  
It demonstrates that transforming the \sde\ to a slow model for the complex amplitude is a delicate process that requires careful treatment of noise integrals in order to form a consistent model of the long term evolution.  
The specific and formal analysis herein needs to be extended to generic oscillatory systems to discover general modelling principles.

\paragraph{Acknowledgements:} I thank Phil Pollett and Yannis
Kevrekidis for useful discussions on aspects of this article.  
This project receives support from the Australian Research Council.

\bibliographystyle{agsm}
\bibliography{ajr,new,bib}

\end{document}